%% file: CS_extension_SS_PapageorgiouNEW.tex
\documentclass[11pt,oneside,reqno]{amsart}

\usepackage[a4paper, total={450pt,675pt}]{geometry}
\usepackage[OT2,T1]{fontenc}
\usepackage[utf8]{inputenc}
\usepackage{amssymb,bm}

\usepackage{float} 

\usepackage[dvipsnames]{xcolor}
\usepackage[
colorlinks=true,
linkcolor=Maroon,
citecolor=JungleGreen,
urlcolor=NavyBlue, linktocpage]{hyperref}

\usepackage[capitalize,nameinlink]{cleveref}

\usepackage{enumitem}
\usepackage{dsfont}
\usepackage{subcaption}
\usepackage{centernot}

\usepackage{yhmath} 

\usepackage{tikz}
\usetikzlibrary{arrows}
\usetikzlibrary{decorations.pathreplacing,angles,quotes}
\usepackage{pgfplots}
\pgfplotsset{compat=1.15}
\usepackage{mathrsfs}
\usetikzlibrary{arrows}

\usepackage{soul}
\setstcolor{red}

\makeatletter
\def\@tocline#1#2#3#4#5#6#7{\relax
	\ifnum #1>\c@tocdepth 
	\else
	\par \addpenalty\@secpenalty\addvspace{#2}%
	\begingroup \hyphenpenalty\@M
	\@ifempty{#4}{%
		\@tempdima\csname r@tocindent\number#1\endcsname\relax
	}{%
		\@tempdima#4\relax
	}%
	\parindent\z@ \leftskip#3\relax
	\advance\leftskip\@tempdima\relax
	\rightskip\@pnumwidth plus4em \parfillskip-\@pnumwidth
	#5\leavevmode\hskip-\@tempdima
	\ifcase #1
	\or\or \hskip 2em \or \hskip 2em \else \hskip 3em \fi%
	#6\nobreak\relax
	\dotfill\hbox to\@pnumwidth{\@tocpagenum{#7}}\par
	\nobreak
	\endgroup
	\fi} 
\makeatother

\usepackage[
backend=bibtex,
style=alphabetic,
maxbibnames=10,
sorting=nyt]{biblatex}

\DeclareLabelalphaTemplate{
	\labelelement{
		\field[final]{shorthand}
		\field{label}
		\field[strwidth=3]{labelname}}
	\labelelement{
		\field[strwidth=2,strside=right]{year}}
}

\DeclareFieldFormat[article,book,incollection,misc]{title}{#1}
\DeclareFieldFormat[inbook,incollection]{booktitle}{#1}
\DeclareFieldFormat[misc]{date}{\textit{preprint} (#1)}
\DeclareFieldFormat{pages}{#1}
\DeclareFieldFormat{doi}{
	\href{http://www.ams.org/mathscinet-getitem?mr=#1}{#1}}

\renewbibmacro{in:}{%
	\ifentrytype{article}{}{\printtext{\bibstring{in}\intitlepunct}}}

\newtheorem{theorem}{Theorem}[section]

\newtheorem{lemma}[theorem]{Lemma}
\newtheorem{proposition}[theorem]{Proposition}
\newtheorem{corollary}[theorem]{Corollary}

\newtheorem{remark}[theorem]{Remark}

\crefname{section}{Sect.}{section}
\numberwithin{equation}{section}

\newcommand*\diff{\mathop{}\!\mathrm{d}}

\DeclareMathOperator{\supp}{supp}

\DeclareMathOperator{\const}{const.}

\allowdisplaybreaks

\begin{document}
	\title[Asymptotic behavior of solutions to the extension problem on  symmetric spaces]{Asymptotic behavior of solutions to the \\ extension problem for the fractional Laplacian\\ on noncompact symmetric spaces}
	
	\author{Effie Papageorgiou}
	
	\begin{abstract}
		This work deals with the extension problem for the fractional Laplacian on Riemannian symmetric spaces 
		$G/K$ of noncompact type and of general rank, which gives rise to a family of convolution operators, including the Poisson operator. More precisely, motivated by Euclidean results for the Poisson semigroup, we study the long-time asymptotic behavior of solutions to the extension problem for $L^1$ initial data.
		In the case of the Laplace-Beltrami operator, we show that if the initial data is bi-$K$-invariant, then the solution to the extension problem behaves asymptotically as the mass times the fundamental solution, but this convergence may break down
		in the non bi-$K$-invariant case. In the second part, we investigate the long-time 
		asymptotic behavior of the extension problem associated 
		with the so-called distinguished Laplacian on $G/K$.
		In this case, we observe phenomena which are similar to the
		Euclidean setting for the Poisson semigroup, such as $L^1$ asymptotic convergence without the assumption of bi-$K$-invariance.
	\end{abstract}
	
	
	\makeatletter
	\@namedef{subjclassname@2020}{\textnormal{2020}
		\it{Mathematics Subject Classification}}
	\makeatother
	
	\maketitle
	\tableofcontents

	\section{Introduction}\label{Section.1 Intro}
	Let $\mathcal{M}$ be a complete, noncompact Riemannian manifold and $\Delta$ be its Laplace-Beltrami operator. It is well understood that the long time behavior of solutions to the heat equation
	\begin{align}\label{S1 heat}
		\begin{cases}
			\partial_{t}u(t,x)&\,
			=\,\Delta u(t,x),
			\qquad\,t>0,\,\,x\in\mathcal{M},\\[5pt]
			u(0,x)&\,=\,u_0(x),
		\end{cases}
	\end{align} is strongly related to the global geometry of $\mathcal{M}$. This applies also to the heat kernel  $h_{t}\left( x,y\right) $, that is, the minimal
	positive fundamental solution of the heat equation or, equivalently, the
	integral kernel of the heat semigroup $\exp \left( t\Delta \right) $ (see
	for instance \cite{Gri2009}).
	
	The connection between the long time behavior of the
	solution $u (t, x)$ of (\ref{S1 heat}) for initial data $u_0\in L^1(\mathcal{M}, \mu)$ (where
	$\mu$ is the Riemannian measure on $\mathcal{M}$) and that of the heat kernel $h_t(x, y)$  has recently been the subject of extensive studies, see for example  \cite{APZ2023, GPZ2022, Vaz2019} or see \cite{AbAl2022, AGMP21, APZ2023, P2023} for variants and related questions.  Denote by  $M=\int_{\mathcal{M}}\diff{\mu(x)}\,u_0(x)$ the mass of the initial data. In the case when $\mathcal{M}=\mathbb{R}^{n}$
	with the Euclidean metric, the heat kernel is given by 
	\begin{equation*}
		h_{t}(x,y)\,=\,(4\pi {t})^{-\frac{n}{2}}e^{-\frac{|x-y|^{2}}{4t}}
	\end{equation*}%
	and the solution to (\ref{S1 heat}) satisfies as $t\rightarrow +\infty $ 
	\begin{equation}
		\Vert u(t,\,.\,)\,-\,M\,h_{t}(\,.\,,x_{0})\Vert _{L^{1}(\mathbb{R}%
			^{n})}\,\longrightarrow \,0  \label{S1 L1 R}
	\end{equation}%
	and 
	\begin{equation}
		t^{\frac{n}{2}}\,\Vert u(t,\,.\,)\,-\,M\,h_{t}(\,.\,,x_{0})\Vert _{L^{\infty }(%
			\mathbb{R}^{n})}\,\longrightarrow \,0.  \label{S1 Linf R}
	\end{equation}%
	By interpolation, a similar convergence holds with respect to any $L^{p}$
	norm when $1<{p}<\infty $: 
	\begin{equation*}
		t^{\frac{n}{2p^{\prime }}}\,\Vert u(t,\,.\,)\,-\,M\,h_{t}(\,.\,,x_0)\Vert
		_{L^{p}(\mathbb{R}^{n})}\,\longrightarrow \,0
	\end{equation*}%
	where $p^{\prime }$ is the H\"{o}lder conjugate of $p$.

	The situation is drastically different in hyperbolic spaces: it was shown by V\'{a}zquez \cite{Vaz2019} that (\ref{S1 L1 R}) fails for a general initial
	function $f\in L^{1}$ but is still true if $%
	f$ is spherically symmetric around $x_{0}.$ Similar results were
	obtained in \cite{APZ2023} in a more general setting of symmetric spaces of
	noncompact type by using tools of harmonic analysis.  In \cite{GPZ2022}, it was shown that (\ref{S1 Linf R}) fails on connected sums $\mathbb{R}^n\#\mathbb{R}^n$, $n\geq 3$.

	The \textit{fractional} Laplacian is the operator $(-\Delta)^{\sigma}$, $\sigma \in(0,1)$, defined as the spectral $\sigma$-th power of the Laplace-Beltrami operator, with $\text{Dom}(-\Delta)\subset\text{Dom}((-\Delta)^{\sigma})$. 
	It is connected to \textit{anomalous} diffusion, which accounts for much of the interest
	in modeling with fractional equations (quasi-geostrophic flows, turbulence and water waves, molecular dynamics,
	and relativistic quantum mechanics of stars). It also has various applications in probability and finance. One can obtain the fractional Laplacian through a Dirichlet-to-Neumann map extension problem,
	introduced by Caffarelli and Silvestre \cite{CaSi2007} on $\mathbb{R}^n$. This extension problem was considered for fractional powers of more general self-adjoint operators in \cite{ST2010}, where also a Poisson formula was given, as well as conditions for the existence of a fundamental solution. On certain ``good'' noncompact Riemannian manifolds $\mathcal{M}$ (e.g. Cartan-Hadamard manifolds or manifolds with non-negative Ricci curvature), the problem was studied in \cite{BanEtAl}.
	More precisely, let $H^{\sigma}(\mathcal{M})$ denote the usual Sobolev space on $\mathcal{M}$. Then, for any given $v_0\in H^{\sigma}(\mathcal{M})$ there exists a unique solution of the extension problem
	\begin{equation}\label{S1 CS problem}
		\Delta v+\frac{(1-2\sigma)}{t}\frac{\partial v}{\partial t}+\frac{\partial^2 v}{\partial t^2}=0, \quad 0 <\sigma < 1,
	\end{equation}
	with $v(0,x)=v_0(x)$ where $t>0$, $x\in \mathcal{M}$, and the fractional Laplacian, can be recovered through
	\begin{equation*}
		(-\Delta)^{\sigma}v_0(x)=-2^{2\sigma-1}\frac{\Gamma(\sigma)}{\Gamma(1-\sigma)}\lim_{t\rightarrow 0^{+}}t^{1-2\sigma}\,\frac{\partial v}{\partial t}(x,t).
	\end{equation*}
	The extension problem has drawn much attention. Since the associated literature is enormous, we shall refer indicatively to \cite{ARBB2022, BanEtAl, BP2022, BrPa2023, FMT13, RT20, ST2010} and the references therein. From a probabilistic point of view, the extension problem corresponds to the property that all symmetric stable processes can be obtained as traces of degenerate Bessel diffusion processes, see \cite{Stinga}.

	Observe that for  $\sigma=1/2$ we get the Poisson semigroup $e^{-t\sqrt{-\Delta}}$. In the Euclidean case $\mathcal{M}=\mathbb{R}^n$, the Poisson kernel is given by 
	\begin{equation}\label{Poisson eucl} \mathcal{Q}_t(x,y)=\frac{\Gamma(\frac{n+1}{2})}{\pi^{\frac{n+1}{2}}}\frac{t}{(t^2+|x-y|^2)^{\frac{n+1}{2}}}, \qquad x\in \mathbb{R}^n, \; t>0.
	\end{equation}
	Then, for $f\in L^1(\mathbb{R}^n)$ as $t\rightarrow +\infty $, it holds \cite{Vaz18F}
	\begin{equation}
		\Vert e^{-t\sqrt{-\Delta}}f\,-\,M\,	\mathcal{Q}_{t}(\,.\,,x_{0})\Vert _{L^{1}(\mathbb{R}%
			^{n})}\,\longrightarrow \,0  \label{S1 L1 R Poisson}
	\end{equation}%
	and 
	\begin{equation}
		t^{n}\,\Vert e^{-t\sqrt{-\Delta}}f\,-\,M\,\mathcal{Q}_{t}(\,.\,,x_{0})\Vert _{L^{\infty }(%
			\mathbb{R}^{n})}\,\longrightarrow \,0.  \label{S1 Linf R Poisson}
	\end{equation}%
	
	Motivated by this, we examine the long time behavior of solutions to the extension problem with $L^1$ initial data on noncompact symmetric spaces $\mathbb{X}=G/K$,  for all $0<\sigma<1$. More precisely, let $Q_t^{\sigma}$ be the fundamental solution to the extension problem \eqref{S1 CS problem}, which will be called \textit{the fractional Poisson kernel} from now on. Then our main result is the following.
	
	\begin{theorem}\label{S1 Main thm 1}
		Let $v_0\in L^1(\mathbb{X})$ be bi-$K$-invariant and consider the solution $v$  to the extension problem (\ref{S1 CS problem}) with initial data $v_0$. Set $M=\int_{\mathbb{X}}v_0.$  Then
		\begin{align}\label{S1 Main thm 1 convergence}
			\|v(t,\,\cdot\,)-M \, Q_{t}^{\sigma}\|_{L^{1}(\mathbb{X})}\,
			\longrightarrow\,0
			\qquad\textnormal{as}\quad\,t\rightarrow+\infty.
		\end{align}
		Moreover, this convergence fails in general without the bi-$K$-invariance assumption.
	\end{theorem}

	To the best of our knowledge, this is the first approach to examine this property on (essentially) negatively curved manifolds. The result is new even for the case of real hyperbolic space.

	\begin{remark}
		If the bi-$K$-invariant initial data is in addition compactly supported, 
		we obtain the better estimate
		\begin{align*}
			\|v(t,\,\cdot\,)-M \, Q_{t}^{\sigma}\|_{L^{1}(\mathbb{X})}\,
			\le\,C\,t^{-\frac{\varepsilon}{2}} 
			\qquad\forall\,t\ge1,
		\end{align*}
		where $C>0$ is a constant and $\varepsilon$ is any
		positive constant such that $\varepsilon<2/(\nu+2\sigma)$,
		see \cref{S3 Sub1}	and \cref{S3 Sub2}. Here, $\nu$ denotes the so-called dimension at infinity 
		of $\mathbb{X}$, see \cref{Section.2 Prelim}.
	\end{remark} 
	
	\begin{remark}
		We also provide the following sup norm (for which no bi-$K$-invariance is needed)
		and $L^{p}$ $(1<p<\infty$) norm estimates:
		\begin{align}
			\|v(t,\,\cdot\,)-\,M\, Q_{t}^{\sigma}\|_{L^{\infty}(\mathbb{X})}\,
			&=\,
			\mathrm{O}\big(t^{-(\frac{\nu}{2}+\frac{1}{2}+\sigma)}e^{-|\rho|t}\big)\label{LinftyX}\\[5pt]
			\|v(t,\,\cdot\,)- \,M\, Q_{t}^{\sigma}\|_{L^{p}(\mathbb{X})}\,
			&=\,\mathrm{o}\big(t^{-\frac{1}{p'}(\frac{\nu}{2}+\frac{1}{2}+\sigma)}
			e^{-\tfrac{|\rho|t}{p'}}\big)
		\end{align}
		as $t\rightarrow +\infty$. 
		Here, $p'$ denotes the dual exponent of $p$, defined by the formula
		$\tfrac{1}{p}+\tfrac{1}{p'}=1$.
		However, the sup norm estimate \eqref{LinftyX} in the present context is relatively weaker compared to
		\eqref{S1 Linf R Poisson} in the Euclidean setting, while the $L^{p}$ norm estimate is 
		similar.
		Here, $\rho$ is the 
		half sum of positive roots with multiplicities, see \cref{Section.2 Prelim}. This is reminiscent of the weak $L^{\infty}$ convergence for the heat equation on $\mathbb{X}$, observed first on three-dimensional real hyperbolic space \cite{Vaz2019} and generalized to arbitrary rank noncompact symmetric spaces \cite{APZ2023}.
	\end{remark}

	Let $S=N(\exp{\mathfrak{a}})=(\exp{\mathfrak{a}})N$ be the solvable group
	occurring in the Iwasawa decomposition $G=N(\exp{\mathfrak{a}})K$. 
	Then $S$ is identifiable, as a manifold, with the symmetric space 
	$\mathbb{X}=G/K$.
	Our second main contribution is to study the asymptotic convergence 
	for solutions to the extension  associated with the 
	so-called distinguished Laplacian
	$\widetilde{\Delta}$ on $S$. 
	In order to state the results, let us introduce some indispensable notation,
	which will be clarified in \cref{Section.2 Prelim}
	and \cref{Section.5 Distinguished}.
	Denote by
	$\varphi_{0}$ the ground spherical function, 
	by $\widetilde{\delta}$ the modular function on $S$, 
	and by $\widetilde{Q}_{t}$
	the fundamental solution to the extension problem
	\begin{align}\label{S1 HE S}
		\widetilde{\Delta} \widetilde{v} - \frac{(1-2\sigma)}{t} \partial_{t} \widetilde{v} - \partial^{2}_{tt} \widetilde{v}=0, \quad 
		\widetilde{v}(\, \cdot \, ,0)=\widetilde{v}_0 ,  \quad t>0.
	\end{align}
	Let $\widetilde{\varphi}_{0}=\widetilde{\delta}^{1/2}\varphi_{0}$ 
	be the modified ground spherical function and denote by 
	$\widetilde{M}=\tfrac{\widetilde{v}_{0}*
		\widetilde{\varphi}_{0}}{\widetilde{\varphi}_{0}}$ the mass function
	on $S$ which generalizes the mass in the Euclidean case (see \cref{Subsection L1 S}).
	Then, we show the following long-time asymptotic convergence results.
	\begin{theorem}\label{S1 Main thm 2}
		Let $\widetilde{v}_{0}$ belong to the class of  continuous and compactly supported functions on 
		$S$. Then, the solution to the extension problem \eqref{S1 HE S} with initial data $\widetilde{v}_{0}$ satisfies
		\begin{align}\label{S1 L1 disting}
			\|\widetilde{v}(t,\,\cdot\,)-
			\widetilde{M}\,\widetilde{Q}_{t}^{\sigma}\|_{L^{1}(S)}\,
			\longrightarrow\,0
		\end{align}
		and
		\begin{align}\label{S1 Linf disting}
			t^{\ell+|\Sigma_{r}^{+}|}
			\|\widetilde{v}(t,\,\cdot\,)-
			\widetilde{M}\,\widetilde{Q}_{t}\|_{L^{\infty}(S)}\,
			\longrightarrow\,0
		\end{align}
		as $t\rightarrow +\infty$. Here $\ell$ denotes the rank of $G/K$ 
		and $\Sigma_{r}^{+}$ the set of positive reduced roots.
		Analogous $L^p$ ($1<p<\infty$) norm estimates 
		follow by interpolation.
	\end{theorem}
	
	\begin{remark}
		Let us comment on \eqref{S1 L1 disting} and \eqref{S1 Linf disting}.
		Firstly, notice that the $L^1$ convergence \eqref{S1 L1 disting} holds 
		without the restriction of  bi-$K$-invariance, in contrast to \cref{S1 Main thm 1}, 
		and the sup norm estimate $\eqref{S1 Linf disting}$ is stronger than
		(\ref{LinftyX}), as in the Euclidean
		setting for the Poisson semigroup.
		Secondly, the mass $\widetilde{M}$ is a bounded function and not necessarily a constant.
		Thirdly, the power $\ell+|\Sigma_{r}^{+}|$ 
		which occurs in time factor, never coincides with the dimension at infinity
		$\nu=\ell+2|\Sigma_{r}^{+}|$ and it is equal to the topological dimension
		$n=\ell+\sum_{\alpha\in\Sigma^{+}}m_{\alpha}$ 
		if and only if the following equivalent conditions hold:
		\begin{itemize}
			\item 
			the root system $\Sigma$ is reduced and all roots have multiplicity
			$m_{\alpha}=1$.
			\item
			$G$ is a normal real form.
		\end{itemize}
	\end{remark}

	This paper is organized as follows. After the present introduction
	in \cref{Section.1 Intro} and preliminaries in \cref{Section.2 Prelim}, we discuss the  extension problem associated with the Laplace-Beltrami operator on symmetric spaces in \cref{Section.3 CS}. In \cref{Section 4 X} we deal with the long-time asymptotic behavior of solutions to the extension problem
	associated with the Laplace-Beltrami operator on symmetric spaces. We first determine the critical region where the fractional Poisson kernel concentrates. Next,
	on the one hand, for continuous compactly supported initial data, we show that both the solution and the fractional Poisson kernel vanish
	asymptotically outside that critical region. On the other hand, inside the critical region we discuss the role of the additional assumption on the bi-$K$-invariance of the initial data. The rest of this section deals with problems for more general initial data in the 
	$L^{p}\,(p\ge1)$ setting. In \cref{Section.5 Distinguished}, we investigate the asymptotic behavior
	of solutions to the extension problem associated with the distinguished Laplacian.
	After specifying the critical region in this context, we study the long-time
	convergence in $L^1$ and in $L^{\infty}$ with compactly supported initial
	data and address some questions associated with other initial data at the end 
	of the paper.

	Throughout this paper, the notation
	$A\lesssim{B}$ between two positive expressions means that 
	there is a constant $C>0$ such that $A\le{C}B$. 
	The notation $A\asymp{B}$ means that $A\lesssim{B}$ and $B\lesssim{A}$. Also, $A(t)\sim B(t)$ means that $A(t)/B(t)\rightarrow 1$ as $t\rightarrow +\infty$.

	\section{Preliminaries}\label{Section.2 Prelim}
	In this section, we review spherical Fourier analysis 
	on Riemannian symmetric spaces of noncompact type. The notation is standard 
	and follows \cite{Hel1978,Hel2000,GaVa1988}. 
	Next we recall bounds and asymptotics of the heat kernel, for which 
	we refer to \cite{AnJi1999,AnOs2003} for more details in this setting.  
	
	\subsection{Noncompact Riemannian symmetric spaces}
	Let $G$ be a semi-simple Lie group, connected, noncompact, with finite center, 
	and $K$ be a maximal compact subgroup of $G$. The homogeneous space 
	$\mathbb{X}=G/K$ is a Riemannian symmetric space of noncompact type.
	Let $\mathfrak{g}=\mathfrak{k}\oplus\mathfrak{p}$ be the Cartan decomposition 
	of the Lie algebra of $G$. The Killing form of $\mathfrak{g}$ induces 
	a $K$-invariant inner product $\langle\,.\,,\,.\,\rangle$ on $\mathfrak{p}$, 
	hence a $G$-invariant Riemannian metric on $G/K$.
	We denote by $d(\,.\,,\,.\,)$ the Riemannian distance on $\mathbb{X}$.
	
	Fix a maximal abelian subspace $\mathfrak{a}$ in $\mathfrak{p}$. 
	The rank of $\mathbb{X}$ is the dimension $\ell$ of $\mathfrak{a}$.
	We identify $\mathfrak{a}$ with its dual $\mathfrak{a}^{*}$ 
	by means of the inner product inherited from $\mathfrak{p}$.
	Let $\Sigma\subset\mathfrak{a}$ be the root system of 
	$(\mathfrak{g},\mathfrak{a})$ and denote by $W$ the Weyl group 
	associated with $\Sigma$. 
	Once a positive Weyl chamber $\mathfrak{a}^{+}\subset\mathfrak{a}$ 
	has been selected, $\Sigma^{+}$ (resp. $\Sigma_{r}^{+}$ 
	or $\Sigma_{s}^{+}$)  denotes the corresponding set of positive roots 
	(resp. positive reduced, i.e., indivisible roots or simple roots).
	Let $n$ be the dimension and $\nu$ be the pseudo-dimension 
	(or dimension at infinity) of $\mathbb{X}$: 
	\begin{align}\label{S2 Dimensions}
		\textstyle
		n\,=\,
		\ell+\sum_{\alpha \in \Sigma^{+}}\,m_{\alpha}
		\qquad\textnormal{and}\qquad
		\nu\,=\,\ell+2|\Sigma_{r}^{+}|
	\end{align}
	where $m_{\alpha}$ denotes the dimension of the positive root subspace
	\begin{align*}
		\mathfrak{g}_{\alpha}\,
		=\,\lbrace{
			X\in\mathfrak{g}\,|\,[H,X]=\langle{\alpha,H}\rangle{X},\quad
			\forall\,H\in\mathfrak{a}
		}\rbrace.
	\end{align*}
	Denote by $\rho\in\mathfrak{a}^{+}$ the half sum of all positive roots 
	$\alpha \in \Sigma^{+}$ counted with their multiplicities $m_{\alpha}$:
	\begin{align*}
		\rho\,=\,
		\frac{1}{2}\,\sum_{\alpha\in\Sigma^{+}} \,m_{\alpha}\,\alpha.
	\end{align*}
	Sometimes we shall use coordinates on $\mathfrak{a}$. When we do, we always refer to the
	coordinates associated to the orthonormal basis $\delta_1, ... , \delta_{\ell-1}, \,\rho/|\rho|$, where $ \delta_1, ... , \delta_{\ell-1},$ is any orthonormal basis of $\rho^{\perp}.$

	Let $\mathfrak{n}$ be the nilpotent Lie subalgebra 
	of $\mathfrak{g}$ associated with $\Sigma^{+}$ 
	and let $N = \exp \mathfrak{n}$ be the corresponding 
	Lie subgroup of $G$. We have the decompositions 
	\begin{align*}
		\begin{cases}
			\,G\,=\,N\,(\exp\mathfrak{a})\,K 
			\qquad&\textnormal{(Iwasawa)}, \\[5pt]
			\,G\,=\,K\,(\exp\overline{\mathfrak{a}^{+}})\,K
			\qquad&\textnormal{(Cartan)}.
		\end{cases}
	\end{align*}
	Denote by $A(x)\in\mathfrak{a}$ and $x^{+}\in\overline{\mathfrak{a}^{+}}$
	the middle components of $x\in{G}$ in these two decompositions, respectively, and by
	$|x|=|x^{+}|$ the distance to the origin.
	In the Cartan decomposition, the Haar measure 
	on $G$ writes
	\begin{align*}
		\int_{G}\diff{x}\,f(x)
		=\,
		|K/\mathbb{M}|\,\,\int_{K}\diff{k_1}\,
		\int_{\mathfrak{a}^{+}}\diff{x^{+}}\,\delta(x^{+})\, 
		\int_{K}\diff{k_2}\,f(k_{1}(\exp x^{+})k_{2})\,,
	\end{align*}
	with density
	\begin{align}\label{S2 estimate of delta}
		\delta(x^{+})\,
		=\,\prod_{\alpha\in\Sigma^{+}}\,
		(\sinh\langle{\alpha,x^{+}}\rangle)^{m_{\alpha}}\,
		\asymp\,
		\prod_{\alpha\in\Sigma^{+}}
		\Big( 
		\frac{\langle\alpha,x^{+}\rangle}
		{1+\langle\alpha,x^{+}\rangle}
		\Big)^{m_{\alpha}}\,
		e^{2\langle\rho,x^{+}\rangle}
		\qquad\forall\,x^{+}\in\overline{\mathfrak{a}^{+}}. 
	\end{align}
	Here $K$ is equipped with its normalized Haar measure, $\mathbb{M}$ denotes the centralizer of $\exp\mathfrak{a}$ in $K$ and the volume 
	of $K/\mathbb{M}$ can be computed explicitly, see \cite[Eq (2.2.4)]{AnJi1999}.

	Finally, let us recall that 
	\begin{align}\label{dist flat}
		|x^{+}-y^{+}|\leq d(xK,yK), \quad |(yx)^{+}-y^{+}|, \;	|(xy)^{+}-y^{+}|\leq d(xK,eK), 
	\end{align}
	see \cite[Lemma 2.1.2]{AnJi1999} or \cite[Lemma 2.1]{MMV2017}.
	
	\subsection{Spherical Fourier analysis} For this subsection, our main references are  \cite[Chap.4]{GaVa1988} and 
	\cite[Chap.IV]{Hel2000}.
	
	For every $\lambda\in\mathfrak{a}$, the spherical function
	$\varphi_{\lambda}$ is a smooth bi-$K$-invariant eigenfunction of all 
	$G$-invariant differential operators on $\mathbb{X}$, in particular of the
	Laplace-Beltrami operator:
	\begin{equation*}
		-\Delta\varphi_{\lambda}(x)\,
		=\,(|\lambda|^{2}+|\rho|^2)\,\varphi_{\lambda}(x).
	\end{equation*}
	It is symmetric in the sense that
	$\varphi_{\lambda}(x^{-1})=\varphi_{-\lambda}(x)$,
	and is given by the integral representation
	\begin{align}\label{S2 Spherical Function}
		\varphi_{\lambda}(x)\, 
		=\,\int_{K}\diff{k}\,e^{\langle{i\lambda+\rho,\,A(kx)}\rangle}.
	\end{align}
	
	All the elementary spherical functions $\varphi_{\lambda}$ 
	with parameter $\lambda\in\mathfrak{a}$ are controlled by the ground spherical
	function $\varphi_{0}$, which satisfies the global estimate
	\begin{align}\label{S2 global estimate phi0}
		\varphi_{0}(\exp{x^{+}})\,
		\asymp\,
		\Big\lbrace \prod_{\alpha\in\Sigma_{r}^{+}} 
		1+\langle\alpha,x^{+}\rangle\Big\rbrace\,
		e^{-\langle\rho, x^{+}\rangle}
		\qquad\forall\,x^{+}\in\overline{\mathfrak{a}^{+}}.
	\end{align}

	Let $\mathcal{S}(K \backslash{G}/K)$ be the Schwartz space of bi-$K$-invariant
	functions on $G$. The spherical Fourier transform (Harish-Chandra transform)
	$\mathcal{H}$ is defined by
	\begin{align}\label{S2 HC transform}
		\mathcal{H}f(\lambda)\,
		=\,\int_{G} \diff{x}\,\varphi_{-\lambda}(x)\,f(x)
		\qquad\forall\,\lambda\in\mathfrak{a},\
		\forall\,f\in\mathcal{S}(K\backslash{G/K}),
	\end{align}
	where $\varphi_{\lambda}\in\mathcal{C}^{\infty}(K\backslash{G/K})$ is the
	spherical function of index $\lambda \in \mathfrak{a}$.
	Denote by $\mathcal{S}(\mathfrak{a})^{W}$ the subspace 
	of $W$-invariant functions in the Schwartz space $\mathcal{S}(\mathfrak{a})$. 
	Then $\mathcal{H}$ is an isomorphism between $\mathcal{S}(K\backslash{G/K})$ 
	and $\mathcal{S}(\mathfrak{a})^{W}$. The inverse spherical Fourier transform 
	is given by
	\begin{align}\label{S2 Inverse formula}
		f(x)\,
		=\,\frac{C_0}{|W|}\,\int_{\mathfrak{a}}\frac{\diff{\lambda}}{|\mathbf{c}(\lambda)|^{2}}\,
		\varphi_{\lambda}(x)\,
		\mathcal{H}f(\lambda) 
		\qquad\forall\,x\in{G},\
		\forall\,f\in\mathcal{S}(\mathfrak{a})^{W},
	\end{align}
	where the constant $C_0=2^{n-\ell}/(2\pi)^{\ell}|K/\mathbb{M}|$ depends only 
	on the geometry of $\mathbb{X}$, and $|\mathbf{c(\lambda)}|^{-2}$ is the 
	so-called Plancherel density, given by an explicit formula by Gindikin-Karpelevič. Finally, if $f$ is a Schwartz function on $\mathbb{X}$, the Helgason-Fourier transform is defined by
	\begin{align}\label{S2 Helgason}
		\widehat{f}(\lambda,k\mathbb{M})\,
		=\,
		\int_{G}\diff{g}\,
		f(gK)\,e^{\langle{-i\lambda+\rho,\,A(k^{-1}g)}\rangle}, 
	\end{align} 
	which, in view of \eqref{S2 Spherical Function}, boils down to the transform \eqref{S2 HC transform}
	when $f$ is bi-$K$-invariant.

	\subsection{Heat kernel on symmetric spaces}
	The heat kernel on $\mathbb{X}$ is a positive bi-$K$-invariant right 
	convolution kernel, i.e., $h_{t}(xK,yK)=h_{t}(y^{-1}x)>0$, 
	which is thus determined by its restriction 
	to the positive Weyl chamber. In fact, it symmetric, i.e. $h_t(x)=h_t(x^{-1})$.
	According to the inversion formula \eqref{S2 Inverse formula} of the spherical Fourier transform, 
	the heat kernel is given by
	\begin{align}\label{S2 heat kernel inv}
		h_{t}(xK)\,
		=\,\frac{C_{0}}{|W|}\,\int_{\mathfrak{a}}\frac{\diff{\lambda}}{|\mathbf{c}(\lambda)|^{2}}\,
		\varphi_{\lambda}(x)\,
		e^{-t(|\lambda|^{2}+|\rho|^{2})}
	\end{align}
	and satisfies the global estimate
	\begin{align}\label{S2 heat kernel}
		h_{t}(\exp{H})\,
		\asymp\,t^{-\frac{n}{2}}\,
		\Big\lbrace{
			\prod_{\alpha\in\Sigma_{r}^{+}}
			(1+t+\langle{\alpha,H}\rangle)^{\frac{m_{\alpha}+m_{2\alpha}}{2}-1}
		}\Big\rbrace\,\varphi_{0}(\exp{H})
		e^{-|\rho|^{2}t-\frac{|H|^{2}}{4t}}
	\end{align}
	for all $t>0$ and $H\in\overline{\mathfrak{a}^{+}}$, 
	see \cite{AnJi1999,AnOs2003}. Recall that $\int_{\mathbb{X}}h_{t}=1$.
	
	Finally, in order to describe more accurately the asymptotic behavior of the ground spherical function and of the heat kernel on certain regions,  let us introduce the following functions: consider 
	\begin{equation}\label{pi function}
		\bm{\pi}(i\lambda)=\prod_{\alpha\in\Sigma_{r}^{+}}
		\langle{\alpha,\lambda}\rangle
	\end{equation}
	and 
	\begin{align*}
		\mathbf{b}(\lambda)\,
		=\,\prod_{\alpha\in\Sigma_{r}^{+}}\,
		\mathbf{b}_{\alpha}
		\Big(
		\frac{\langle\alpha,\lambda\rangle}{\langle\alpha,\alpha\rangle}
		\Big)
	\end{align*}
	where 
	\begin{align*}
		\mathbf{b}_{\alpha}(z)\,
		=\,|\alpha|^{2}\,
		\tfrac{\Gamma(\frac{\langle{\alpha,\rho}\rangle}
			{\langle{\alpha,\alpha}\rangle}
			+\frac{1}{2} m_{\alpha})}
		{\Gamma(\frac{\langle{\alpha,\rho}\rangle}
			{\langle{\alpha,\alpha}\rangle})}\,
		\tfrac{\Gamma(\frac{1}{2}
			\frac{\langle{\alpha,\rho}\rangle}
			{\langle{\alpha,\alpha}\rangle} 
			+\frac{1}{4} m _{\alpha} 
			+ \frac{1}{2} m_{2\alpha})}
		{\Gamma(\frac{1}{2}
			\frac{\langle{\alpha,\rho}\rangle}
			{\langle {\alpha,\alpha}\rangle} 
			+ \frac{1}{4} m_{\alpha})}\,
		\tfrac{\Gamma(iz+1)}
		{\Gamma(iz+ \frac{1}{2}m_{\alpha})}\,
		\tfrac{\Gamma(\frac{i}{2}z 
			+ \frac{1}{4} m_{\alpha})}
		{\Gamma(\frac{i}{2}z 
			+\frac{1}{4} m_{\alpha} 
			+ \frac{1}{2} m_{2\alpha})}.
	\end{align*}
	The function $\mathbf{b}(-\lambda)^{-1}$ is holomorphic for 
	$\lambda\in\mathfrak{a}+i\overline{\mathfrak{a}^{+}}$ and positive for 
	$\lambda\in i\overline{\mathfrak{a}^{+}}$. 
	We recall that it has the following
	behavior
	\begin{align}\label{bfunction}
		|\mathbf{b}(-\lambda)|^{-1}\,
		\asymp\,
		\prod_{\alpha\in\Sigma_{r}^{+}}\,
		(1+|\langle{\alpha,\lambda}\rangle|)^{\frac{m_{\alpha}+m_{2\alpha}}{2}-1}
	\end{align}
	and that its derivatives can be estimated by
	\begin{align}\label{bfunction derivative}
		p(\tfrac{\partial}{\partial\lambda})
		\mathbf{b}(-\lambda)^{-1}\,
		=\,\mathrm{O}\big(|\mathbf{b}(-\lambda)|^{-1}\big),
	\end{align}
	where $p(\tfrac{\partial}{\partial\lambda})$ is any differential 
	polynomial, \cite[pp.1041-42]{AnJi1999}.
	
	We are now ready to describe the asymptotic behavior of the ground spherical function away from the walls. More precisely, as $\mu(H)=\min_{\alpha\in\Sigma^{+}}\langle{\alpha,H}\rangle\rightarrow\infty$, we have 
	\begin{align}\label{S2 phi0 far}
		\varphi_{0}(\exp{H})\,
		\sim\,C_{1}\,\bm{\pi}(H)\,e^{-\langle{\rho,H}\rangle}.
	\end{align}
	Here, $C_{1}=\bm{\pi}(\widetilde{\rho})^{-1}\mathbf{b}(0)$ and $\widetilde{\rho}=\frac{1}{2}\sum_{\alpha\in \Sigma_{r}^{+}}\alpha$, see \cite[Proposition 2.2.12(ii)]{AnJi1999}.
	
	As for the heat kernel, we have the 
	following
	asymptotics \cite[Theorem 5.1.1]{AnJi1999}:
	\begin{align}\label{AJ heat asymp}
		h_t(\exp{H})\,
		\sim\,C_{2}\,t^{-\frac{\nu}{2}}\,
		\mathbf{b}\big(-i\tfrac{H}{2t}\big)^{-1}\,
		\varphi_{0}(\exp{H})\,e^{-|\rho|^{2}t-\frac{|H|^{2}}{4t}}
	\end{align}
	as $t\rightarrow\infty$, provided
	$\mu(H)\rightarrow\infty$ or $|H|=\textrm{O}(t)$.
	Here $C_{2}=C_{0}2^{-|\Sigma_{r}^{+}|}\pi^{\frac{\ell}{2}}
	\bm{\pi}(\widetilde{\rho})\mathbf{b}(0)^{-1}$.

	\section{The fractional Laplacian and the extension problem} \label{Section.3 CS}
	This section deals with the notion of the fractional Laplacian and the extension problem which gives rise to a family of operators, containing the Poisson operator.

	In recent years there has been intensive research on various kinds of fractional order operators. Being nonlocal objects, local PDE techniques to treat nonlinear problems for the fractional operators do not apply. To overcome this difficulty, in the Euclidean case, Caffarelli and Silvestre \cite{CaSi2007} studied the extension problem associated with the Laplacian and realized the fractional power as the map taking Dirichlet data to Neumann data. In \cite{ST2010} Stinga and Torrea related the extension problem for the fractional Laplacian to the heat semigroup. On certain classes of noncompact manifolds, which include symmetric spaces of noncompact type, the extension problem has been studied by Banica, Gonz{\'a}lez and S{\'a}ez \cite{BanEtAl}. Interestingly, in the noncompact setting one needs to have a precise control of the behavior of the metric at infinity and geometry plays a crucial role. 
	
	From now on, we strictly work on symmetric spaces of noncompact type $\mathbb{X}=G/K$. To begin with, using the spectral theorem, one can define fractional powers of the Laplacian via the heat semigroup,
	$$(-\Delta)^{\sigma}f(x)=\int_{0}^{\infty}\frac{\diff t}{t^{1+\sigma}}\,(e^{t\Delta}f(x)-f(x)) \quad \text{ in } L^2(\mathbb{X}), \; f\in \text{Dom}(-\Delta). $$
	see \cite[(5), p.260]{Y}, \cite{BanEtAl}, or \cite{ST2010}.

	Then, by \cite[Theorems 1.1 and 2.1]{ST2010} (see also \cite[Theorem 1.1]{BanEtAl}) the relation between the fractional Laplacian and the extension problem (\ref{S1 CS problem}) is the following.
	\begin{theorem} \cite{ST2010}. Let $\sigma\in (0,1)$. Then for  $f\in \text{Dom}((-\Delta)^{\sigma})$, a solution to the extension problem
		\begin{align*}
			\Delta v+\frac{(1-2\sigma)}{t}\frac{\partial v}{\partial t}+\frac{\partial^2 v}{\partial t^2}=0, \quad  
			v(0,x)\,=\,f(x),\quad  t>0,\,\; x\in \mathbb{X},
		\end{align*}
		is given by
		\begin{equation*}
			v(t,x)=v(t, g K)=f\ast Q_t^{\sigma}(g K)=\int_G \diff{y}\, Q_t^{\sigma}(y^{-1}g)\,f(y), \quad g\in G,
		\end{equation*}
		where 
		\begin{equation}\label{Q kernel}
			Q_t^{\sigma}(g)=\frac{t^{2\sigma}}{2^{2\sigma}\Gamma(\sigma)}\int_{0}^{+\infty}\frac{\diff u}{u^{1+\sigma}} \,h_u(g)\,e^{-\frac{t^2}{4u}},
		\end{equation}
		Moreover, the fractional Laplacian  on $\mathbb{X}$ can be recovered through
		\begin{equation*}(-\Delta)^{\sigma}f(x)=-2^{2\sigma-1}\frac{\Gamma(\sigma)}{\Gamma(1-\sigma)}\lim_{t\rightarrow 0^{+}}t^{1-2\sigma}\frac{\partial v}{\partial t}(x,t).
		\end{equation*}
	\end{theorem}
	
	It is worth mentioning that the existence of an integral (in fact, a right convolution) kernel on symmetric spaces for the extension problem  follows from \cite[pp.2099-2101]{ST2010} (see also \cite[Theorem 3.2]{BanEtAl}), since the heat kernel 
	for $t\in (0,1)$ and all $x\in \mathbb{X}$ satisfies
	\begin{equation}\label{L2heat}
		\|h_t(x,\, \cdot \, )\|_2 \lesssim t^{-\frac{n}{4}}, \quad \|\partial_t h_t(x, \, \cdot \, )\|_{2}\lesssim t^{-\frac{n}{4}-1}.
	\end{equation}
	This follows by an explicit computation using the heat kernel estimates \eqref{S2 heat kernel}, along with \eqref{S2 global estimate phi0} and \eqref{S2 estimate of delta}, while for the time derivative, one may use the pointwise estimate \cite[(3.1)]{A92}.
	
	Observe that due to the subordination \eqref{Q kernel} to the heat kernel, $Q_t^{\sigma}$ is a positive, bi-$K$-invariant and symmetric  (in the sense that $Q_t^{\sigma}(g)=Q_t^{\sigma}(g^{-1})$ for all $g\in G$) function on $G$.

	We next recall some large-time upper and lower bounds  for the kernel $Q_{t}^{\sigma}$ proved in \cite{BP2022}.
	\begin{theorem}\cite[Theorem 3.2]{BP2022}
		The fractional Poisson kernel $Q_t^\sigma$ on $\mathbb{X}$, $0<\sigma<1$,  satisfies the following upper and lower bounds 
		\begin{equation}\label{P-bound}
			Q_t^{\sigma} (\exp H)\asymp  \frac{t^{2\sigma}}{4^\sigma\,\Gamma(\sigma)} \left(\sqrt{t^2+|H|^2}\right)^{-\frac{\ell}{2}-\frac{1}{2}-\sigma-|\Sigma_r^+|} \varphi_0(\exp H) ~e^{-|\rho|\sqrt{t^2+|H|^2}},
		\end{equation}
		if $|H|^2+t^2 \geq 1$.
	\end{theorem}

	For large time sup norm estimates, we have the following result.
	\begin{proposition}\label{sup norm est}
		For $t>1$, it holds
		$$\|Q_t^{\sigma}\|_{L^{\infty}(\mathbb{X})}\asymp t^{\sigma-\frac{\ell}{2}-\frac{1}{2}-|\Sigma_r^+|}e^{-|\rho|t}.$$
	\end{proposition}
	\begin{proof}
		The lower bound follows immediately by the fact that $\|Q_t^{\sigma}\|_{L^{\infty}(\mathbb{X})}\geq Q_t^{\sigma}(eK)$, \eqref{P-bound} and the fact that $\varphi_0(eK)=1$. For the upper bound, we use that 
		$$\left(\sqrt{t^2+|H|^2}\right)^{-\frac{\ell}{2}-\frac{1}{2}-\sigma-|\Sigma_r^+|}\leq t^{-\frac{\ell}{2}-\frac{1}{2}-\sigma-|\Sigma_r^+|}, \quad e^{-|\rho|\sqrt{t^2+|H|^2}}\leq e^{-|\rho|t},  $$
		and the fact that $\varphi_{0}(\exp H)\lesssim 1$, for all $H\in\overline{\mathfrak{a}^{+}}$.
	\end{proof}
	
	\subsection{Large time behavior of the fractional Poisson kernel}	Recall first that $\int_{G}\diff{g}\,h_{t}(g)=1,$  $\forall t>0$. This implies that 
	$$\int_{G}\diff{g}\,Q_{t}^{\sigma}(g)=1, \quad \forall t>0,\; \forall \sigma \in (0,1),$$ 
	by the subordination formula \eqref{Q kernel}, the definition of the Gamma function and a Fubini argument. Motivated by this, let us introduce the notion of the \textit{critical region} for the kernel $Q_t^{\sigma}$.
	
	\begin{proposition}\label{prop critical region} Let $0<\varepsilon<1$. Consider in $\mathfrak{a}$ the annulus
		\begin{equation} 
			t^{2-\varepsilon}\leq |H|\leq t^{2+\varepsilon}
		\end{equation}
		and the solid cone $\Gamma(t)$ with angle $$\gamma(t)=t^{-\frac{\varepsilon}{2}}$$ around the $\rho$-axis, and denote by  $\Omega_t$ their intersection. Then, the critical region for the fractional Poisson kernel is $K(\exp\Omega_{t})K$, in the sense that 
		$$\int_{G \smallsetminus K(\exp\Omega_{t})K}\diff{g}\,Q_{t}^{\sigma}(g)\longrightarrow 0, \quad \text{as } t\rightarrow +\infty.$$
	\end{proposition}
	
	\begin{figure}
		\centering
		\input{Concentration1Lnew.tex}
		\caption{Flat part $\Omega_t$ of critical region}
		\label{fig concentration1}
	\end{figure}
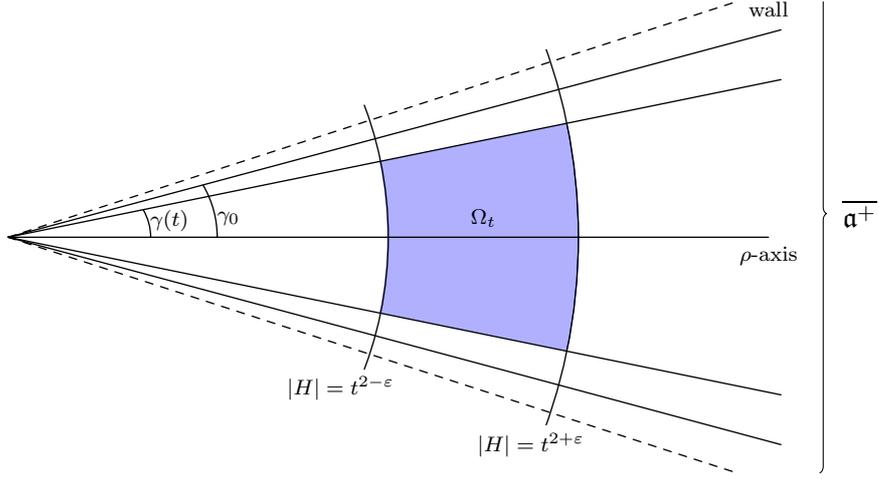
	
	\begin{proof}
		Let the rank $\ell$ be greater or equal to $2$ (the rank case one case is simpler, thus omitted). Let $0\leq a<b$. Using the bounds \eqref{P-bound} and the fact that 
		\begin{equation*}
			\varphi_0(\exp H)\lesssim (1+|H|)^{|\Sigma_r^{+}|}e^{-\langle \rho, H \rangle},
		\end{equation*} 
		we have, by the Cartan decomposition and \eqref{S2 estimate of delta},
		\begin{align}\label{integral in H}
			& \int_{a<|x|<b} \diff x \, Q_t^\sigma (x) \nonumber \\
			& \lesssim t^{2\sigma}\int_{\{a\leq |H|\leq b\}\cap\overline{\mathfrak{a}^+}}\diff H\, \left(\sqrt{t^2+|H|^2}\right)^{-\frac{\ell}{2}-\frac{1}{2}-\sigma-|\Sigma_r^{+}|}\,(1+|H|)^{|\Sigma_r^{+}|}\,e^{\langle \rho, H \rangle}\,e^{-|\rho| \sqrt{t^2+|H|^2}}.
		\end{align}
		
		Take $t$ large enough so that $\Gamma(t)$ is contained inside a small cone $\Gamma_0$ with fixed angle $\gamma_0$ around the $\rho$-axis
		and consider the regions 
		\begin{align*}
			R_1&=\{x\in G: \, |x|< t^{2-\varepsilon}\}, \\
			R_2&=\{x\in G: \, |x|> t^{2-\varepsilon},\; x\notin K(\exp\Gamma_0)K\},\\
			R_3&=\{x\in G: \, t^{2-\varepsilon}\leq |x|\leq t^{2+\varepsilon},\; x\in  K(\exp\Gamma_0)K \smallsetminus K(\exp\Gamma(t))K  \},\\
			R_4&=\{x\in G: \, |x|> t^{2+\varepsilon}, \; x\in  K(\exp\Gamma_0)K \}.
		\end{align*}
		
		First of all, we have
		\begin{equation*}
			\int_{ R_1}\diff x\, Q_t^\sigma(x) \lesssim t^{2\sigma} e^{-\frac{|\rho|}{3}t^{\varepsilon}}\int_{0}^{t^{2-\varepsilon}} \diff r\, (1+ r)^{|\Sigma_r^+|}~r^{\ell-1}  \lesssim t^{-N\varepsilon} \quad \forall N>0, 
		\end{equation*}
		using \eqref{integral in H} and that if $|H|<t^{2-\varepsilon}$ then for $t$ large enough,
		\begin{equation*}
			e^{\langle \rho, H \rangle}\,e^{-|\rho| \sqrt{t^2+|H|^2}} \leq \exp\left\{-|\rho|\frac{t^2}{\sqrt{t^2+|H|^2}+|H|}\right\}\leq \exp\left\{-\frac{|\rho|}{3}t^{\varepsilon}\right\}.
		\end{equation*}
		
		Next, observe that if $H\notin \Gamma_0$, it  holds
		\begin{equation*}
			e^{\langle \rho, H \rangle}\,e^{-|\rho| \sqrt{t^2+|H|^2}}\leq e^{-|\rho||H|(1-\cos\gamma_0)},
		\end{equation*} 
		which yields from \eqref{integral in H} that 
		\begin{align*}
			\int_{R_2} \diff x\,Q_t^{\sigma} (x)&\lesssim t^{2\sigma} \int_{\{H\notin \Gamma_0:\;|H|> t^{2-\varepsilon}\}} \diff H\,(1+|H|)^{|\Sigma_r^+|} e^{-|\rho||H|(1-\cos \gamma_0)}\\
			&\lesssim  t^{2\sigma} e^{-\frac{|\rho|}{2}\,t^{2-\varepsilon}(1-\cos \gamma_0)} \lesssim t^{-N\varepsilon} \quad \forall N>0.
		\end{align*}

		We next pass to the region $R_3$. Recall first the trivial inequality 
		\begin{equation*}
			\sin \theta \geq \frac{2}{\pi}\theta, \quad \theta \in [0,\pi/2].
		\end{equation*}
		Then, for some positive constant $c=c(|\rho|)>0$, we have  that
		\begin{align*}
			e^{\langle \rho, H \rangle}\,e^{-|\rho| \sqrt{t^2+|H|^2}} &\leq  e^{-|\rho||H|(1-\cos \gamma(t))}\\
			&= e^{-2|\rho||H|\sin^2 (\gamma(t)/2)}\\
			&\leq e^{-2|\rho||H| \gamma(t)^2/\pi^2}\\
			&\leq e^{-c\,t^{2-2\varepsilon}},
		\end{align*}
		since $|H|\geq t^{2-\varepsilon}$ and  $\gamma(t)=t^{-\frac{\varepsilon}{2}}$. Since $0<\varepsilon<1$, passing to polar coordinates,  we get by \eqref{integral in H} 
		\begin{align*}
			\int_{R_3}\diff x\, Q_t^\sigma (x) &\lesssim t^{2\sigma}\int_{t^{2-\varepsilon}}^{t^{2+\varepsilon}} \diff r\, \left(\sqrt{t^2+r^2}\right)^{-\frac{\ell}{2}-\frac{1}{2}-\sigma-|\Sigma_r^{+}|}\,(1+r)^{|\Sigma_r^{+}|}\,r^{\ell-1} e^{-c\,t^{2-2\varepsilon}} \\
			&\lesssim t^{-N\varepsilon} \quad \forall N>0.
		\end{align*}
		
		To treat the integral in the remaining region $R_4$, in view of \eqref{integral in H}, let us write in polar coordinates:
		\begin{align}
			\int_{R_4}\diff x\, Q_t^\sigma(x)   
			&\lesssim t^{2\sigma} \int_{t^{2+\varepsilon}}^{+\infty} \diff r\, \left(\sqrt{t^2+r^2}\right)^{-\frac{\ell}{2}-\frac{1}{2}-\sigma-|\Sigma_r^{+}|}\,(1+r)^{|\Sigma_r^{+}|}\,r^{\ell-1} \nonumber\\ 
			& \times \int_{0}^{\gamma_0} \diff\gamma\, e^{-|\rho|\frac{t^2}{\sqrt{t^2+r^2}+r\cos \gamma}}e^{-|\rho|\frac{r^2\sin^2\gamma}{\sqrt{t^2+r^2}+r\cos \gamma}}\,\sin^{\ell-2}\gamma. \label{polar}
		\end{align}
		
		Observe that		
		\begin{align}
			\int_{0}^{\gamma_0} \diff\gamma\, e^{-|\rho|\frac{r^2\sin^2\gamma}{\sqrt{t^2+r^2}+r\cos \gamma}}\,\sin^{\ell-2}\gamma
			&\leq  \int_{0}^{\gamma_0} \diff\gamma\, e^{-|\rho| \frac{r^2 \gamma^2 4/\pi^2}{\sqrt{t^2+r^2}+r}}\,\gamma^{\ell-2} \notag\\
			&\lesssim  \left(\frac{r^2}{\sqrt{t^2+r^2}+r} \right)^{\frac{1-\ell}{2}}.	\label{polar contrib} 
		\end{align}

		Thus, in $R_4$, where $r^2+t^2\asymp r^2$, \eqref{polar} and \eqref{polar contrib} yield
		\begin{align*}
			\int_{R_4} \diff x\,Q_t^\sigma (x) &\lesssim_\sigma t^{2\sigma} \int_{ t^{2+\varepsilon}}^{+\infty}\diff r\, r^{-\frac{\ell}{2}-\frac{1}{2}-\sigma-|\Sigma_r^{+}|}(1+r)^{|\Sigma_r^{+}|}~r^{\ell-1}r^{\frac{1-\ell}{2}}\\
			& \lesssim t^{2\sigma} \int_{ t^{2+\varepsilon}}^{+\infty} \diff r \, r^{-\sigma-1}\\
			&\lesssim t^{-\sigma \varepsilon}.
		\end{align*}
		This completes the proof.
	\end{proof}
	\begin{remark} The corresponding critical region for the Poisson kernel $\mathcal{Q}_t$ in the Euclidean case would be $B(0,t^{1+\varepsilon})\smallsetminus B(0,t^{1-\varepsilon})$, as one can easily check using \eqref{Poisson eucl}. On the other hand, the heat kernel $h_t$ on a Riemannian symmetric space of the noncompact	type  is asymptotically concentrated along the ($K$-orbit) of the $\rho$-axis and an annulus centered at the origin, however moving to
		infinity with finite speed $2|\rho|$, \cite{AnSe1992}.
	\end{remark}
	
	We now obtain  precise long-time asymptotics of the kernel $Q_t^{\sigma}$ which are crucial for our proof, by a slightly more general result.
	
	\begin{theorem} \label{asyp-P}
		Let $\sigma \in (0, 1)$. Then, as $t+|H|\rightarrow+\infty$, we have
		\begin{align}
			Q_t^\sigma(\exp H) &\sim  C(\sigma)\,
			t^{2\sigma}\, \bm{b}\left(-i |\rho|\frac{ H}{\sqrt{t^2+|H|^2}}\right)^{-1} \left(\sqrt{t^2+|H|^2}\right)^{-\frac{\ell}{2}-\sigma-|\Sigma_r^+|-\frac{1}{2}}\times \notag\\
			&\times \varphi_0(\exp H)\, e^{-|\rho| \sqrt{t^2+|H|^2}}~,
		\end{align}
		where the constant is 
		$$C(\sigma)=\frac{1}{4^\sigma \Gamma(\sigma)}~C_0~2^{\frac{\ell}{2}+\sigma+\frac{1}{2}}~\pi^{\frac{\ell}{2}+\frac{1}{2}} \bm{\pi}(\widetilde \rho)\,\bm{b}(0)^{-1}|\rho|^{\frac{\ell}{2}+\sigma+|\Sigma_r^+|-\frac{1}{2}},$$
		with $\widetilde{\rho}=\frac{1}{2}\sum_{\alpha \in \Sigma_r^{+}}\alpha$ and $C_0=2^{n-\ell}/(2\pi)^{\ell}|K/\mathbb{M}|$.
	\end{theorem}	
	
	\begin{proof}
		The proof follows arguments for the asymptotics of the Poisson kernel ($\sigma=1/2$) in \cite[Section 5]{AnJi1999}. 
		
		Consider a constant $\kappa>4$. In view of the subordination formula \eqref{Q kernel}, let us split 
		\begin{align*}
			Q_t^\sigma(x) &= \frac{t^{2\sigma}}{4^\sigma \Gamma(\sigma)}\int_{0}^{+\infty} \frac{\diff u}{u^{1+\sigma}} \,h_u(x)\,e^{-\frac{t^2}{4u}},\\
			&=\frac{t^{2\sigma}}{4^\sigma \Gamma(\sigma)}\{J_1+J_2+J_3\},
		\end{align*}
		where the quantities $J_1, J_2$ and $J_3$ are defined by the integration  over the intervals $\big[0, \kappa^{-1} b \big), \big[\kappa^{-1} b, \kappa b\big)$ and $\big[\kappa b, \infty\big)$ respectively. Here,  $b=\frac{\sqrt{t^2+|x|^2}}{2|\rho|}$.
		
		We claim that the main contribution comes from the middle integral $J_2$. Indeed,  for the first integral $J_1$, we get that for some $\delta>0$ and some constants $d_1, d_2>0$, we have
		\begin{align*}
			J_1 \lesssim (1+|x|)^{d_2-|\Sigma_r^{+}|}(t^2+|x|^2 )^{-\sigma-d_1}\,\varphi_0(x)\, e^{-(|\rho|+\delta)\sqrt{t^2+|x|^2}}  ,
		\end{align*}
		\cite[p.19]{BP2022}. For the third integral $J_3$, we get again by \cite[p.18]{BP2022} that
		\begin{align*}
			J_3\lesssim \left(\sqrt{t^2+|x|^2}\right)^{-\frac{\ell}{2}-|\Sigma_r^+|-\frac{1}{2}} \varphi_0(x)\,e^{-(|\rho|+\eta)\sqrt{t^2+|x|^2}},
		\end{align*}
		where $\eta= |\rho| \kappa/4 -|\rho|>0$.
		
		We now consider $J_2$. Define
		\begin{align}\label{S2 heat kernel critical region}
			h(t,H)= t^{\frac{\ell}{2}+|\Sigma_{r}^{+}|}\,
			\mathbf{b}\big(-i\tfrac{H}{2t}\big)\,
			\varphi_{0}(\exp{H})^{-1}\,e^{|\rho|^{2}t+\frac{|H|^{2}}{4t}}\, h_t(\exp{H}).
		\end{align}
		Then,  by a change of variables and \eqref{S2 heat kernel critical region}, we have
		\begin{align*}
			J_2&= \int_{\kappa^{-1}\frac{ \sqrt{t^2+|x|^2}}{2|\rho|}}^{\kappa\frac{ \sqrt{t^2+|x|^2}}{2|\rho|}} \frac{\diff u}{u^{1+\sigma}}\, h_u(x)\,e^{-\frac{t^2}{4u}}\\
			&= \left(\frac{\sqrt{t^2+|H|^2}}{2|\rho|}\right)^{-\sigma} \int_{\kappa^{-1}}^{\kappa} \frac{\diff u}{u^{1+\sigma}}\, h_{u\sqrt{t^2+|H|^2}/{2|\rho|}}(\exp H)~e^{-t^2|\rho|/\left(2u\sqrt{t^2+|H|^2}\right)}= \\
			&=\varphi_0(\exp H)~\left(\frac{\sqrt{t^2+|H|^2}}{2|\rho|}\right)^{-\frac{\ell}{2}-\sigma-|\Sigma_{r}^{+}|} \\
			& \times \int_{\kappa^{-1}}^{\kappa} \frac{\diff u}{u^{\frac{\ell}{2}+|\Sigma_r^{+}|-1+\sigma}} \,  \textbf{b} \left(-i \frac{|\rho| H}{u\sqrt{t^2+|H|^2}}\right)^{-1} e^{-\frac{|\rho|}{2}u\sqrt{t^2+|H|^2}-\frac{|\rho|}{2u}\sqrt{t^2+|H|^2}}\, h\left( u\frac{\sqrt{t^2+|H|^2}}{2|\rho|}, H\right)\\
			&	=\varphi_0(\exp H)~\left(\frac{\sqrt{t^2+|H|^2}}{2|\rho|}\right)^{-\frac{\ell}{2}-\sigma-|\Sigma_{r}^{+}|}\textbf{b} \left(-i \frac{|\rho| H}{\sqrt{t^2+|H|^2}}\right)^{-1} \\
			& \times \int_{\kappa^{-1}}^{\kappa} \frac{\diff u}{u^{\frac{\ell}{2}+|\Sigma_r^{+}|-1+\sigma}} \, \frac{\textbf{b} \left(-i \frac{|\rho| H}{u\sqrt{t^2+|H|^2}}\right)^{-1}}{\textbf{b} \left(-i \frac{|\rho| H}{\sqrt{t^2+|H|^2}}\right)^{-1}} e^{-|\rho|\sqrt{t^2+|H|^2}(\frac{u+u^{-1}}{2})}\, h\left( u\frac{\sqrt{t^2+|H|^2}}{2|\rho|}, H\right).
		\end{align*}
		By the Laplace method we get that the last integral tends to 
		$$C_2\,\sqrt{\frac{2\pi}{|\rho|}}\,(t^2+|H|^2)^{-\frac{1}{2}}e^{-|\rho|\sqrt{t^2+|H|^2}}, $$
		due the facts that
		\begin{itemize}
			\item[(i)] $\textbf{b} \left(-i \frac{|\rho| H}{u\sqrt{t^2+|H|^2}}\right)^{-1}$ is bounded above and below, uniformly in $u$ and $t$, $H$;
			\item[(ii)] $\textbf{b} \left(-i \frac{|\rho| H}{u\sqrt{t^2+|H|^2}}\right)^{-1}\sim \textbf{b} \left(-i \frac{|\rho|H}{\sqrt{t^2+|H|^2}}\right)^{-1}$ as $u \rightarrow 1$, uniformly in $t$ and $H$;
			\item[(iii)]	$h\left( u\frac{\sqrt{t^2+|H|^2}}{2|\rho|}, H\right)\longrightarrow C_2$ as $t+|H|\rightarrow+\infty$, uniformly in $u$, by the asymptotics  in \eqref{AJ heat asymp} and contradiction (see \cite[pp.1085-1086]{AnJi1999}).
			
		\end{itemize}
		For the exact value of the constant $C_2$, we refer again to \eqref{AJ heat asymp}.
		
		Since $J_1, J_3$ are very small compared to $J_2$ for $t$ large, substituting the value of $C_2$ we finally get the claimed asymptotics.
	\end{proof}
	
	\subsection{Asymptotics in the critical region}\label{S3 Sub1}
	
	In this subsection we prove asymptotics for some quantities that will be later used in the proof. 
	
	The first two lemmas describe the effect of a small translation on the critical region.
	
	\begin{lemma}\label{lemma asymptotics}
		For all $x$ in the critical region $K(\exp\Omega_{t})K$, and for all $y\in G$ bounded, 
		the following asymptotic behaviors hold as $t\rightarrow +\infty$: 
		\begin{itemize}
			\item[(i)] 
			$\frac{|(y^{-1}x)^{+}|}{|x^{+}|}$ and $\frac{|x^{+}|}{|(y^{-1}x)^{+}|}$ are both equal to 
			$1+\textnormal{O}\big(t^{-2+\varepsilon}\big)$.
			
			\item[(ii)]          
			$\frac{x^{+}}{|x^{+}|}$ and $\frac{(y^{-1}x)^{+}}{|(y^{-1}x)^{+}|}$
			are both equal to
			$\frac{\rho}{|\rho|}+\textnormal{O}\big( t^{-\frac{\varepsilon}{2}}\big)$. 
			
			\item[(iii)]
			For every $\alpha\in\Sigma^{+}$,
			$\frac{\langle{\alpha,(y^{-1}x)^{+}}\rangle}{
				\langle{\alpha,x^{+}}\rangle}
			=1+\textnormal{O}\big( t^{-\frac{\varepsilon}{2}}\big)$.
			
			\item[(iv)] 
			$d(xK,eK)-d(xK,yK)
			=\langle{\frac{\rho}{|\rho|},A(k^{-1}y)}\rangle
			+\textnormal{O}\big( t^{-\frac{\varepsilon}{2}}\big)$. 	Here,  $k$ is the left component of $x$ in the Cartan decomposition and $A(k^{-1}y)$ is the flat middle component of $k^{-1}y$ in the Iwasawa decomposition. 
		\end{itemize}
	\end{lemma}
	
	\begin{proof}
		Assume that $t^{2-\varepsilon}\leq d(xK, eK)\leq t^{2+\varepsilon}$ and $d(yK, eK)\leq \xi$, which implies by the triangle inequality that $\frac{1}{2}\,t^{2-\varepsilon}\leq d(xK, yK)\leq 2\,t^{2+\varepsilon}$, for $t$ large enough.
		
		We deduce first $(i)$ by using
		\begin{align*}
			\frac{|(y^{-1}x)^{+}|}{|x^{+}|}\,
			= \frac{d(xK,yK)}{d(xK, eK)}
			=\,1+\textrm{O}\big(t^{-2+\varepsilon}\big).  
		\end{align*}
		The second assertion follows similarly.
		
		Next, for (ii), since the angle of $x^{+}$ with the $\rho$-axis is $\textrm{O}(t^{-\frac{\varepsilon}{2}})$, we first have
		\begin{align*}
			\left|\frac{x^{+}}{|x^{+}|}-\frac{\rho}{|\rho|} \right|^2
			=2\left(1-\frac{\langle \rho, x^{+}\rangle}{|\rho||x^{+}|}\right)=\textrm{O}(t^{-\varepsilon}).
		\end{align*}
		For the second asymptotics in (ii), we work similarly, observing that since $(y^{-1}x)^{+}=x^{+}+\textrm{O}(1)$, we have
		\begin{align}\label{star}
			\langle \frac{\rho}{|\rho|}, \frac{(y^{-1}x)^{+}}{|(y^{-1}x)^{+}|} \rangle & =\frac{|x^{+}|}{|(y^{-1}x)^{+}|} \langle  \frac{\rho}{|\rho|},\frac{x^{+}}{|x^{+}|}  \rangle +\textrm{O}(|x|^{-1}) \notag \\
			&= \langle  \frac{\rho}{|\rho|},\frac{x^{+}}{|x^{+}|}  \rangle +\textrm{O}(t^{-2+\varepsilon})\\
			&=1+\textrm{O}(t^{-\varepsilon}), \notag
		\end{align}
		using (i) and that $\cos(\widehat{x^{+}, \rho})=1+\textrm{O}(t^{-\varepsilon}).$
		
		Let us next deduce (iii) from $(i)$ and $(ii)$. 
		For every positive root $\alpha$,
		\begin{align*}
			\frac{\langle{\alpha,(y^{-1}x)^{+}}\rangle}{\langle{\alpha,x^{+}}\rangle}\,
			&=\,
			\frac{\langle{\alpha,\frac{(y^{-1}x)^{+}}{|(y^{-1}x)^{+}|}}\rangle}{
				\langle{\alpha,\frac{x^{+}}{|x^{+}|}}\rangle}\,
			\frac{|(y^{-1}x)^{+}|}{|x^{+}|}\notag\\[5pt]
			&=\, 
			\frac{\langle{\alpha,\frac{\rho}{|\rho|}}\rangle
				+\textrm{O}\big( t^{-\frac{\varepsilon}{2}}\big)}{
				\langle{\alpha,\frac{\rho}{|\rho|}}\rangle
				+\textrm{O}\big( t^{-\frac{\varepsilon}{2}}\big)}\,
			\Big\lbrace{1+\textrm{O}\big(t^{-2+\varepsilon}\big)}\Big\rbrace\,
			=\,1+\textrm{O}( t^{-\frac{\varepsilon}{2}}).
		\end{align*}
		
		It remains to prove (iv). For that, we follow \cite[Lemma 3.8]{APZ2023}. Let $x=k(\exp{x^{+}})k'$
		in the Cartan decomposition and consider the Iwasawa decomposition
		$k^{-1}y=n(k^{-1}y)(\exp{A(k^{-1}y)})k''$ for some $k''\in{K}$. Then
		\begin{align*}\label{S3 distance decomposition in lemma}
			d(xK,yK)\,
			&=\,d\big(k(\exp{x^{+}})K,kn(k^{-1}y)(\exp{A(k^{-1}y)})K\big)
			\notag\\[5pt]
			&=\,d\big(\exp{(-x^{+})}[n(k^{-1}y)]^{-1}(\exp{x^{+}})K,
			\exp{(A(k^{-1}y)}-x^{+})K\big).
		\end{align*}
		and we write
		\begin{align*}
			d(xK,eK)-d(xK,yK)\,
			&=\,\overbrace{\vphantom{\Big|}
				d(xK,eK)-d\big(\exp{(A(k^{-1}y)}-x^{+})K,eK\big)}^{I}\\
			&+\,\underbrace{\vphantom{\Big|}
				d\big(\exp{(A(k^{-1}y)}-x^{+})K,eK\big)-d(xK,yK)}_{II}.
		\end{align*}
		On the one hand, $|II|$ tends exponentially fast to $0$, see \cite{APZ2023}.  On the other hand, we have
		\begin{align*}
			I\,
			=\,|x^{+}|-|A(k^{-1}y)-x^{+}|\,
			&=\,\frac{2\langle{x^{+},A(k^{-1}y)}\rangle-|A(k^{-1}y)|^{2}}{
				|x^{+}|+|A(k^{-1}y)-x^{+}|}\\[5pt]
			&=\,\big\langle{\tfrac{x^{+}}{|x^{+}|},A(k^{-1}y)}\big\rangle\,
			+\textrm{O}\big(\tfrac{1}{|x^{+}|}\big)\\[5pt]
			&=\,\big\langle{\tfrac{\rho}{|\rho|},A(k^{-1}y)}\big\rangle\,+\textrm{O}\big( t^{-\frac{\varepsilon}{2}}\big)
		\end{align*}
		by using $(ii)$, the fact that $\lbrace{A(k^{-1}y)\,|\,k\in{K}}\rbrace$ 
		is a compact subset of $\mathfrak{a}$ and that $\textrm{O}\big(\tfrac{1}{|x^{+}|}\big)=\textrm{O}\big(t^{-2+\varepsilon}\big)$. This concludes the proof.
	\end{proof}

	\begin{lemma}\label{inner prod translation}
		Let $x\in K(\exp \Omega_{t})K $ and let $y$ be bounded. Then
		\begin{equation*}
			\langle \rho, x^{+}\rangle -\langle \rho,(y^{-1}x)^{+}\rangle=|\rho||x^{+}|-|\rho||(y^{-1}x)^{+}|+\textrm{O}( t^{-\frac{\varepsilon}{2}}).
		\end{equation*}
	\end{lemma}	
	
	\begin{proof}
		Since the rank one case is trivial, let us consider $\ell\geq 2$. Observe first that the claim follows by 
		\begin{align}\label{cos trans}
			\cos(\widehat{(y^{-1}x)^{+}, \rho})=  \cos(\widehat{x^{+}, \rho})+\textrm{O}(t^{-\frac{\varepsilon}{2}}|x^{+}|^{-1}).
		\end{align}
		Indeed, by \eqref{cos trans} and taking into account that $\cos(\widehat{x^{+}, \rho})=1+\textrm{O}(t^{-\varepsilon})$, we get
		\begin{align*}
			\langle \rho, x^{+}\rangle -\langle \rho,(y^{-1}x)^{+}\rangle&= |\rho||x^{+}|\cos(\widehat{x^{+}, \rho})- |\rho||(y^{-1}x)^{+}|\cos(\widehat{(y^{-1}x)^{+}, \rho})\\
			&=|\rho|\,|x^{+}|\cos(\widehat{x^{+}, \rho})-|\rho|\,|(y^{-1}x)^{+}|\,(\cos(\widehat{x^{+}, \rho})+\textrm{O}(t^{-\frac{\varepsilon}{2}}|x^{+}|^{-1})) \\
			&=|\rho||x^{+}|-|\rho||(y^{-1}x)^{+}| +\textrm{O}(t^{-\frac{\varepsilon}{2}}).
		\end{align*}
		
		Therefore, it remains to prove \eqref{cos trans}. Observe that by \eqref{star}, we have
		\begin{equation*}
			\sin^2\left(\frac{(\widehat{(y^{-1}x)^{+}, \rho}) }{2}\right)=\sin^2\left(\frac{\widehat{(x^{+}, \rho)}}{2}\right)+\textrm{O}(t^{-2+\varepsilon}),
		\end{equation*} 
		whence $(\widehat{(y^{-1}x)^{+}, \rho})=\textrm{O}(t^{-\frac{\varepsilon}{2}})$.
		Next, recall the coordinates on $\mathfrak{a}$ with respect to the basis $\delta_1, ..., \delta_{\ell-1}, \rho/|\rho|$ introduced in Section \ref{Section.2 Prelim}, and write
		$$x^{+}=\left(\xi,\, \xi_{\ell}\right), \qquad (y^{-1}x)^{+}=\left(\zeta, \, \zeta_{\ell}\right).$$
		Since $\langle x^{+}, \rho \rangle=\xi_{\ell}\,|\rho|$, we get
		\begin{equation*}
			\xi_{\ell}=|x^{+}|\cos(\widehat{x^{+}, \rho}), \qquad |\xi|=|x^{+}|\sin(\widehat{x^{+}, \rho}). 
		\end{equation*} 
		Similarly, 
		\begin{equation*}
			\zeta_{\ell}=|(y^{-1}x)^{+}|\cos(\widehat{(y^{-1}x)^{+}, \rho}), \qquad |\zeta|=|(y^{-1}x)^{+}|\sin(\widehat{(y^{-1}x)^{+}, \rho}).
		\end{equation*}
		Therefore
		\begin{equation}\label{xi coord}	
			\frac{|\xi|}{\xi_{\ell}}=\tan(\widehat{x^{+}, \rho})=\textrm{O}(t^{-\frac{\varepsilon}{2}}), \qquad	|x^{+}|\asymp \xi_{\ell}
		\end{equation}
		and 
		\begin{equation}\label{zeta coord}	
			\frac{|\zeta|}{\zeta_{\ell}}=\tan(\widehat{(y^{-1}x)^{+}, \rho})=\textrm{O}(t^{-\frac{\varepsilon}{2}}), \qquad	|(y^{-1}x)^{+}|\asymp \zeta_{\ell}.	
		\end{equation}
		Thus, we have
		\begin{align}\label{cos frac}
			\left| \cos(\widehat{(y^{-1}x)^{+}, \rho})- \cos(\widehat{x^{+}, \rho})  \right||x^{+}|&=\frac{\left|\zeta_{\ell}\,|x^{+}|-\xi_{\ell}\,|(y^{-1}x)^{+}|\right|}{|(y^{-1}x)^{+}|} \notag \\
			&\asymp \frac{\left| \zeta_{\ell}\sqrt{|\xi|^2+\xi_{\ell}^2}-\xi_{\ell}\sqrt{|\zeta|^2+\zeta_{\ell}^2}\right|}{\zeta_{\ell}} \notag \\
			&\asymp \frac{\left| \zeta_{\ell}^2\,|\xi|^2 - \xi_{\ell}^2\,|\zeta|^2\right|}{\zeta_{\ell}^2\,\xi_{\ell}\left( \sqrt{\left( \frac{|\xi|}{\xi_\ell}\right)^2+1} +\sqrt{\left( \frac{|\zeta|}{\zeta_\ell}\right)^2+1} \right)} \notag \\
			&\asymp \frac{\left| \zeta_{\ell}\,|\xi| - \xi_{\ell}\,|\zeta|\right| \left| \zeta_{\ell}\,|\xi| + \xi_{\ell}\,|\zeta|\right|}{\zeta_{\ell}^2\,\xi_{\ell}},
		\end{align}	
		due to \eqref{xi coord} and \eqref{zeta coord}. The fact that $(y^{-1}x)^{+}=x^{+}+\textrm{O}(1)$ implies $\zeta_{\ell}=\xi_{\ell}+\textrm{O}(1)$ and $|\zeta|=|\xi|+\textrm{O}(1)$, therefore
		\begin{equation}\label{cos frac asymp}
			\frac{ \zeta_{\ell}\,|\xi| - \xi_{\ell}\,|\zeta| }{\xi_{\ell}}=\textrm{O}(1), \qquad \frac{ \zeta_{\ell}\,|\xi| + \xi_{\ell}\,|\zeta|}{\zeta_{\ell}^2}=\textrm{O}\left(\frac{|\zeta|}{\zeta_{\ell}}\right)=\textrm{O}(t^{-\frac{\varepsilon}{2}}).
		\end{equation}
		Altogether, we conclude  by \eqref{cos frac} and \eqref{cos frac asymp} that 
		$$ \cos(\widehat{(y^{-1}x)^{+}, \rho})- \cos(\widehat{x^{+}, \rho})= \textrm{O}(t^{-\frac{\varepsilon}{2}}|x^{+}|^{-1}).  $$
	\end{proof}

	The next lemma is the heart of the proof.
	
	\begin{lemma}\label{quotient} 
		Assume that $x=k\exp(x^{+})k'$ is in the critical region $K(\exp\Omega_{t})K$ and that $y$ is bounded. Then,  
		$$\frac{Q_t^{\sigma}(xK,yK)}{Q_t^{\sigma}(xK,eK)}=e^{\langle 2\rho, A(k^{-1}y)\rangle}+\textrm{O}(t^{-\frac{\varepsilon}{2}}),$$
		for $0<\varepsilon<2/(\nu+2\sigma)$. Here, $\nu=\ell+2|\Sigma_{r}^{+}|$ is the dimension at infinity.
	\end{lemma}
	
	\begin{proof}
		By \eqref{asyp-P}, we get
		\begin{align*} 
			\frac{Q_t^{\sigma}(y^{-1}x)}{Q_t^{\sigma}(x)}\sim
			\frac{\textbf{b}\left(-i|\rho|\frac{(y^{-1}x)^{+}}{\sqrt{t^2+|(y^{-1}x)^{+}|^2}}\right)^{-1} \sqrt{t^2+| (y^{-1}x)^{+}|^2}^{\,-\sigma-\frac{\ell}{2}-|\Sigma_r^{+}|-\frac{1}{2}}\,e^{-|\rho| \sqrt{t^2+|(y^{-1}x)^{+}|^2}} \,\varphi_0(y^{-1}x)}{\textbf{b}\left(-i|\rho|\frac{x^{+}}{\sqrt{t^2+|x^{+}|^2}}\right)^{-1}\sqrt{t^2+| x^{+}|^2}^{\,-\sigma-\frac{\ell}{2}-|\Sigma_r^{+}|-\frac{1}{2}}\,e^{-|\rho| \sqrt{t^2+|x^{+}|^2}}\,\varphi_0(x)}.
		\end{align*}
		
		Our aim is to show that for $x$ inside the critical region and $y$ bounded, the following asymptotics hold;
		
		\begin{itemize}	
			\item[(i)] $\frac{\textbf{b}\left(-i|\rho|\frac{(y^{-1}x)^{+}}{\sqrt{t^2+|(y^{-1}x)^{+}|^2}}\right)^{-1}}{
				\textbf{b}\left(-i|\rho|\frac{x^{+}}{\sqrt{t^2+|x^{+}|^2}}\right)^{-1}}=1 + \textnormal{O}\big(t^{-\frac{\varepsilon}{2}}\big)$;
			
			\item[(ii)] $\left( \frac{t^2+|x^{+}|^2}{t^2+|(y^{-1}x)^{+}|^2}\right)^{k}=1+\textnormal{O}\left(t^{-2+\varepsilon(\nu+2\sigma)}\right)$, where $k=\frac{1}{2}\left(\sigma+\frac{\nu}{2}+\frac{1}{2} \right)$;
			
			\item[(iii)] $\exp\left\{-|\rho| \left(\sqrt{t^2+|(y^{-1}x)^{+}|^2}-\sqrt{t^2+|x^{+}|^2}\right)\right\}=\langle{\rho,A(k^{-1}y)}\rangle+\textrm{O}( t^{-2 +2\varepsilon} )$;
			
			\item[(iv)] $\frac{\varphi_{0}(y^{-1}x)}{\varphi_{0}(x)}=\langle{\rho,A(k^{-1}y)}\rangle+\textnormal{O}\big( t^{-\frac{\varepsilon}{2}}\big)$;
		\end{itemize}		
		whence the claim for $Q_t^{\sigma}(y^{-1}x)/Q_t^{\sigma}(x)$ follows,  choosing $\varepsilon$ small enough. 
		
		We start the proof of (i)-(iv) by some preliminary observations. Write 
		$$r=|x^{+}|=d(xK,eK), \qquad s=|(y^{-1}x)^{+}|=d(xK, yK),$$ and let $d(yK,eK)<\xi$, for some $\xi>0$. Then, for $x$ inside the critical region, and $t$ large enough we have $$t^{2-\varepsilon}\leq r\leq t^{2+\varepsilon}, \qquad \frac{1}{2}\,t^{2-\varepsilon}\leq s \leq 2\,t^{2+\varepsilon}, \qquad |r-s|\leq \xi.$$ Also, 
		we have $t^2+r^2\asymp r^2$ and $t^2+s^2\asymp s^2$. Finally, in the proof of  Lemma \ref{inner prod translation} it was shown that the angle of $(y^{-1}x)^{+}$ with the $\rho$-axis is $\textrm{O}(t^{-\frac{\varepsilon}{2}})$.
		
		\textit{Proof of (i).} Observe first that owing to \eqref{bfunction}, for all $g\in G$ such that $t^2+|g^{+}|^2\asymp |g^{+}|^2$, we have
		\begin{equation}\label{bfunc comp}
			\textbf{b}\left(-i|\rho|\frac{g^{+}}{\sqrt{t^2+|g|^2}}\right)^{-1}\asymp 1.
		\end{equation} Therefore, for $x\in K(\exp \Omega_t)K$ and $y$ bounded, by the mean value theorem we get
		\begin{align*} 
			\left|	\textbf{b}\left(-i\frac{|\rho|x^{+}}{\sqrt{t^2+|x^{+}|^2}}\right)^{-1}-\textbf{b}\left(-i\frac{|\rho|(y^{-1}x)^{+}}{\sqrt{t^2+|(y^{-1}x)^{+}|^2}}\right)^{-1} \right| \lesssim \left|\frac{x^{+}}{\sqrt{t^2+|x^{+}|^2}}-\frac{(y^{-1}x)^{+}}{\sqrt{t^2+|(y^{-1}x)^{+}|^2}}\right|,
		\end{align*} 
		where we have used the derivative bound \eqref{bfunction derivative} and \eqref{bfunc comp}. Next, owing to Lemma \ref{lemma asymptotics}(ii), we obtain
		\begin{align*}
			\frac{x^{+}}{\sqrt{t^2+|x^{+}|^2}}&=\frac{x^{+}}{|x^{+}|}\left(1+\frac{t^2}{|x^{+}|^2}\right)^{-1/2}=\left(\frac{\rho}{|\rho|}+\textnormal{O}\big( t^{-\frac{\varepsilon}{2}}\big)\right) \left( 1+\textrm{O}(t^{-2+2\varepsilon})\right)\\
			&=\frac{\rho}{|\rho|}+\textnormal{O}\big( t^{-\frac{\varepsilon}{2}}),
		\end{align*}  
		since $0<\varepsilon<2/(\nu+2\sigma)<2/3$. Working likewise for $\frac{(y^{-1}x)^{+}}{\sqrt{t^2+|(y^{-1}x)^{+}|^2}}$, the claim follows using  \eqref{bfunc comp}. 
		
		\textit{Proof of (ii).} We use a similar mean value argument applied to $(t^2+(.)^2)^{k}$, $k>1$, so that for some $r_0$ between $r$ and $s$ we have
		\begin{equation}\label{MVTk}
			\left|\frac{(t^2+r^2)^{k}}{(t^2+s^2)^{k}}-1 \right|\lesssim \frac{r_0\,(t^2+r_0^2)^{k-1}}{(t^2+s^2)^k}.
		\end{equation}
		Given that $r_0\lesssim t^{2+\varepsilon}$, $s\gtrsim t^{2-\varepsilon}$ and $k=\frac{1}{2}\left(\sigma+\frac{\nu}{2}+\frac{1}{2} \right)$ we get the desired result.

		\textit{Proof of (iii).} We first claim that 
		\begin{equation}\label{exp2}
			\frac{r+s}{\sqrt{t^2+r^2}+\sqrt{t^2+s^2}}=1+O\left(t^{ -2 +2\varepsilon}\right).
		\end{equation}
		Indeed, consider the function $f(\tau)=\sqrt{\tau^2+r^2}+\sqrt{\tau^2+s^2}$, $\tau \geq 0$, and observe that the left hand side of \eqref{exp2} is equal to $f(0)/f(t)$. Then, the mean value theorem for $f$ in $[0,t]$  together with the fact that  $$f'(\tau)\lesssim  \frac{\tau}{r}+\frac{\tau}{s}\lesssim t^{-1+\varepsilon}, \qquad f(\tau)\gtrsim  t^{2-\varepsilon}, \qquad \forall\tau\in[0,t],$$   yield the claimed asymptotics \eqref{exp2}. 
		Finally, in Lemma \ref{lemma asymptotics}(iv) it was shown that
		\begin{equation}\label{exp1}
			r-s=d(xK,eK)-d(xK,yK)
			=\langle \frac{\rho}{|\rho|}, A(k^{-1}y) \rangle
			+\textnormal{O}\big(t^{-2+\varepsilon}\big).
		\end{equation}
		Therefore, by \eqref{exp2} and \eqref{exp1}, we get
		\begin{align*}
			\exp\left\{ -|\rho|\left(\sqrt{t^2+s^2}-\sqrt{t^2+r^2}\right)\right\}
			&=\exp\left\{|\rho|(r-s)\frac{r+s}{\sqrt{t^2+r^2}+\sqrt{t^2+s^2}}\right\} \\
			&=e^{ \langle \rho, A(k^{-1}y)\rangle +\textrm{O}(t^{ -2 +2\varepsilon})}\\
			&=e^{ \langle \rho, A(k^{-1}y)\rangle}+\textrm{O}(t^{ -2 +2\varepsilon }),
		\end{align*}
		which proves (iii).
		
		\textit{Proof of (iv).} Since the angles of both $x^{+}$ and $(y^{-1}x)^{+}$ with the $\rho$-axis are $\textrm{O}(t^{-\frac{\varepsilon}{2}})$ we may use the ground spherical asymptotics \eqref{S2 phi0 far}.
		On the one hand, by Lemma \ref{lemma asymptotics}(iii), we have
		\begin{align}\label{pi asymp}
			\frac{\bm{\pi}((y^{-1}x)^{+})}{\bm{\pi}(x^{+})}\,
			=\,\prod_{\alpha\in\Sigma_{r}^{+}}
			\frac{\langle{\alpha,(y^{-1}x)^{+}}\rangle}{\langle{\alpha,x^{+}}\rangle}\,
			=\,1+\textrm{O}\big( t^{-\frac{\varepsilon}{2}}\big). 
		\end{align}
		On the other hand, using  Lemma \ref{inner prod translation} and Lemma \ref{lemma asymptotics}(iv),  we have 
		\begin{align*}	e^{\langle \rho, x^{+}\rangle -\langle \rho,(y^{-1}x)^{+}\rangle }&=e^{ \langle \rho, A(k^{-1}y)\rangle +\textrm{O}(t^{-2+\varepsilon})}=e^{ \langle \rho, A(k^{-1}y)\rangle}+\textrm{O}(t^{-2+\varepsilon}),
		\end{align*}
		with which the proof of (iv) is complete.
		
		Altogether, we have 
		$$\frac{Q_t^{\sigma}(xK,yK)}{Q_t^{\sigma}(xK,eK)}=e^{ \langle 2\rho, A(k^{-1}y)\rangle}+\textrm{O}( t^{-\frac{\varepsilon}{2}}).$$
	\end{proof}	
	
	\section{Asymptotic convergence associated with 
		the extension problem for the Laplace-Beltrami operator}\label{Section 4 X}
	We first consider continuous compactly supported initial data $v_0$. We work separately outside and inside the critical region: we will show that $$\|v_0 \ast Q_{t}^{\sigma}-M\, Q_t^{\sigma}\|_{L^1(G \smallsetminus K(\exp\Omega_t)K)}\rightarrow 0$$  but inside $K(\exp\Omega_t)K$,  unless $v_0$ is bi-$K$-invariant, the convergence to the fundamental solution may break down.

	\subsection{Estimates outside the critical region}\label{S3 Sub2}
	In this subsection, we show that the solution $v(t,x)$ to the extension problem
	vanishes asymptotically 
	in $L^{1}(G\smallsetminus{K(\exp\Omega_{t})K})$ as $t\rightarrow +\infty$. Then the desired convergence follows by the triangle inequality.
	
	\begin{lemma}\label{critical region translation} Let $x\in G\smallsetminus K(\exp\Omega_{t})K$ and $y\in K(\exp B(0,\xi))K$. Denote by $\Gamma''(t)$ the solid cone around the $\rho$-axis of angle $\frac{1}{2}\, t^{-\frac{\varepsilon}{2}}$. Consider in $\mathfrak{a}$ the set $$\Omega_t''=\left(B(0, 2\,t^{2+\varepsilon})\smallsetminus B\left(0, \frac{1}{2}\,t^{2-\varepsilon}\right) \right)\cap \Gamma''(t).$$
		Then, $$y^{-1}x \in G\smallsetminus K(\exp\Omega_t'')K.$$
	\end{lemma}

	\begin{proof} Let $x\in G\smallsetminus K(\exp\Omega_{t})K$ and $|y|<\xi$. Recall that by \eqref{dist flat}
		\begin{align*}
			|(y^{-1}x)^{+}-x^{+}|\,
			\le\,d(yK,eK)\,
			=\,|y|\,<\,\xi,
		\end{align*} 
		which implies that
		\begin{align*}
			\begin{cases}
				|(y^{-1}x)^{+}|\,\le\,|x^{+}|+\xi\,
				<\,t^{2-\varepsilon}+\xi\,
				<\,2\,t^{2-\varepsilon},\\
				|(y^{-1}x)^{+}|\,\ge\,|x^{+}|-\xi\,
				>\,t^{2+\varepsilon}-\xi >\,\frac{1}{2}\,t^{2+\varepsilon}.
			\end{cases}
		\end{align*}
		for $t$ large enough. In other words,
		\begin{align*}
			x\,\in\,G\smallsetminus{K(\exp \{B(0, t^{2+\varepsilon})\smallsetminus B(0, t^{2-\varepsilon})\})K}
			\;\Longrightarrow\;
			y^{-1}x\,\in\,G\smallsetminus{K(\exp \{B(0, \frac{1}{2}t^{2+\varepsilon})\smallsetminus B(0, 2\,t^{2-\varepsilon})\})K}.
		\end{align*}
		
		We finally turn to the angles. Write $\phi=(\widehat{x^{+}, \rho})$ and $\omega=(\widehat{(y^{-1}x)^{+}, \rho})$, and observe that by \eqref{star}, we have
		\begin{equation}\label{phi omega}
			\sin^2\left(\frac{\phi}{2}\right)=\sin^2\left(\frac{\omega}{2}\right)+\textrm{O}(t^{-2+\varepsilon}).
		\end{equation}
		Using that 
		$$\sin^2\left(\frac{\phi}{2}\right)\geq \frac{1}{\pi^2}\,\phi^2\geq \frac{1}{\pi^2}\,t^{-\varepsilon}, \qquad \sin^2\left(\frac{\omega}{2}\right)\leq \frac{1}{4}\,\omega^2,$$
		we get that $\omega\geq \frac{1}{2}t^{-\frac{\varepsilon}{2}}$, for $t$ large enough. This completes the proof.
	\end{proof}
	
	\begin{proposition}
		The solution to the extension problem  satisfies
		\begin{align}\label{S3 solution estimate outside}
			\|v(t,\,\cdot\,)\|_{L^{1}(G\smallsetminus{K(\exp\Omega_{t})K})}\,
			\lesssim\,t^{-\sigma\varepsilon}
		\end{align}
		for $t>0$ large enough.
	\end{proposition}

	\begin{proof}
		Let $\xi>0$ be a constant such that the compact support of $v_0$
		belongs to $K(\exp{B(0,\xi)})K$. Then,
		\begin{align*}
			\int_{G\,\smallsetminus\,K(\exp\Omega_{t})K}\,\diff{x}\,|v(t,x)|\, 
			&\lesssim\,
			\int_{K(\exp{B(0,\xi)})K}\diff{y}\,|v_0(y)|\,
			\int_{G\,\smallsetminus\,K(\exp\Omega_{t})K}\diff{x}\,Q_{t}^{\sigma}(y^{-1}x)\\
			&\lesssim\,
			\int_{K(\exp{B(0,\xi)})K}\diff{y}\,|v_0(y)|\,
			\int_{G\,\smallsetminus\,K(\exp\Omega_{t}'')K}\diff{z}\,Q_{t}^{\sigma}(z)
		\end{align*}
		where  $\Omega_t''\subseteq\mathfrak{a}$ is the region described in Lemma \ref{critical region translation}.  
		Thus, working as in Proposition \ref{prop critical region}, one can show that the right-hand side of 
		the inequality above is  $\textrm{O}(t^{-\sigma\varepsilon})$. In conclusion,
		\begin{align*}
			\int_{G\,\smallsetminus\,K(\exp\Omega_{t})K}\diff{x}\,|v(t,x)|\,
			\lesssim\,
			t^{-\sigma\varepsilon}.
		\end{align*}
	\end{proof}

	\subsection{Long-time behavior inside the critical region}
	
	Let now $x\in K(\exp\Omega_t)K$. By Lemma \ref{quotient}, the right-$K$-invariance of $A(k^{-1}.)$ and $v_0$, and the definition \eqref{S2 Helgason} of the Helgason-Fourier transform we have that
	\begin{align}
		v_0 \ast Q_{t}^{\sigma}(x)-M\,Q_{t}^{\sigma}(x)&= \int_{G}\diff{y} \,(Q_{t}^{\sigma}(y^{-1}x)-Q_{t}^{\sigma}(x))v_0(y)\notag \\
		&=Q_{t}^{\sigma}(x)\int_{G}\diff{y} \,\left(\frac{Q_{t}^{\sigma}(y^{-1}x)}{Q_{t}^{\sigma}(x)}-1\right)v_0(y) \notag \\
		&=Q_{t}^{\sigma}(x)\left\{\int_{G}\diff{y} \,\left(e^{ \langle2 \rho, A(k^{-1}y) 
			\rangle}-1+\textrm{O}( t^{-\frac{\varepsilon}{2}})\right)v_0(y)\right\} \notag \\
		&=
		Q_t^{\sigma}(x)\,
		\big(
		\widehat{v_0}(i\rho,k\mathbb{M})\,
		-\,\widehat{v_0}(-i\rho,k\mathbb{M})\,
		+\,\textrm{O}\big( t^{-\frac{\varepsilon}{2}}\big)
		\big).
		\label{S3 new rmk}
	\end{align}
	Notice that $\widehat{v_0}(\pm\,i\rho,k\mathbb{M})=\mathcal{H}v_0(\pm\,i\rho)=M$ when $v_0$ is bi-$K$-invariant. Then we deduce the desired convergence by integrating \eqref{S3 new rmk} over the critical region:
	\begin{equation}\label{convergence bi-K}
		\int_{K(\exp\Omega_t)K}\,
		\diff{x}\,|v_0 \ast Q_{t}^{\sigma}(x)\,-\,M\, Q_{t}^{\sigma} (x)|\, =\textrm{O}( t^{-\frac{\varepsilon}{2}}).
	\end{equation}
	
	On the other hand, using again the Cartan decomposition  we have
	\begin{align*}
		\int_{K(\exp\Omega_t)K}\,
		\diff{x}\,|v_0 \ast Q_{t}^{\sigma}(x)\,-\,M\, Q_{t}^{\sigma} (x)|\, 
		&\longrightarrow\,
		\int_{K}\diff{k}\,
		\Big|
		\int_{G}\diff{y}\,v_0(y)
		\big(
		e^{ \langle 2\rho, A(k^{-1}y) 
			\rangle}\, \,-\,1
		\big) \Big|
	\end{align*}
	as $t\rightarrow +\infty$. The last integral is not constantly zero when $v_0$ is not bi-$K$-invariant. For example, consider $v_0$ to be a Dirac measure supported on some point $yK$ other than the origin, thus for $y\notin K$. In other words,  the solution now coincides with $Q_t(.,yK)$ and the mass is equal to $1$. In this case, however, the last integral is equal to $\int_{K}\diff{k}\,
	\Big|e^{\langle 2\rho , A(k^{-1}y)\rangle}\, \,-\,1
	\,\Big|$, thus does not vanish identically.
	
	\subsection{Long-time convergence for general bi-$K$-invariant data}
	In this subsection,	using the results of the previous two ones and a standard density argument, we prove Theorem \ref{S1 Main thm 1} for the whole class of $L^{1}(\mathbb{X}$) functions that are bi-$K$-invariant. The argument is identical to that of \cite[Subsection 3.3]{APZ2023} but we include it for the reader's convenience.

	\begin{proof}[Proof of  Theorem \ref{S1 Main thm 1}]
		Let $\varepsilon>0$, $v_0\in L^1(K\backslash{G}/K)$ and $V_{0}\in\mathcal{C}_{c}^{\infty}(K\backslash{G}/K)$ 
		be such that $\|v_{0}-V_{0}\|_{L^{1}(\mathbb{X})}<\tfrac{\varepsilon}{3}$.
		Denote by $M=\int_{G}v_0$ and  $M_{V}=\int_{G}V_{0}$ the 
		masses of $v_0$ and $V_{0}$ respectively, then
		\begin{align*}
			|M-M_{V}|\,
			\le\,\|v_{0}-V_{0}\|_{L^{1}(\mathbb{X})}\,
			<\,\tfrac{\varepsilon}{3}.
		\end{align*}
		Let $V(t,x)=V_{0}*Q_{t}^{\sigma}(x)$ be the solution to
		the extension problem with initial data $V_{0}$. We deduce from \eqref{convergence bi-K}, \eqref{S3 solution estimate outside} and
		Proposition \ref{P-bound} that, there exists $T>0$ such that
		\begin{align*}
			\|V(t,\,\cdot\,)-M_{V}\,Q_{t}^{\sigma}\|_{L^{1}(\mathbb{X})}
			&\le\,\|V(t,\,\cdot\,)-M_{V}\,Q_{t}^{\sigma}\|_{L^{1}(K(\exp\Omega_{t})K)}\,
			+\,\|V(t,\,\cdot\,)\|_{L^{1}(G\smallsetminus{K(\exp\Omega_{t})K})}\\[5pt]
			&+\,|M_{V}|\,
			\|Q_{t}^{\sigma}\|_{L^{1}(G\smallsetminus{K(\exp\Omega_{t})K})}\\[5pt]
			&<\,\tfrac{\varepsilon}{3}
		\end{align*}
		for all $t\ge{T}$. In conclusion,
		\begin{align*}
			\|v(t,\,\cdot\,)-M\,Q_{t}^{\sigma}\|_{L^{1}(\mathbb{X})}
			&\le\overbrace{\vphantom{\Big|}
				\|v(t,\,\cdot\,)-V(t,\,\cdot\,)\|_{L^{1}(\mathbb{X})}
			}^{\le\,\|v_{0}-V_{0}\|_{L^{1}}\,\|Q_{t}^{\sigma}\|_{L^{1}}}
			+\overbrace{\vphantom{\Big|}
				\|M_{V}\,Q_{t}^{\sigma}-M\,Q_{t}^{\sigma}\|_{L^{1}(\mathbb{X})}
			}^{\le\,|M-M_{V}|\,\|Q_{t}^{\sigma}\|_{L^{1}}}\\[5pt]
			&+\|V(t,\,\cdot\,)-M_{V}\,Q_{t}^{\sigma}\|_{L^{1}(\mathbb{X})}
			\\[5pt]
			&<\,\tfrac{\varepsilon}{3}+\tfrac{\varepsilon}{3}+\tfrac{\varepsilon}{3}\,
			=\varepsilon
		\end{align*}
		for all $\varepsilon>0$ and $t$ large enough.
	\end{proof}
	
	\vspace{5pt}
	Let us turn to the long-time convergence in $L^{p}(\mathbb{X})$ with $p>1$. We first deal with the case $p=\infty$ and conclude for all $1<p<\infty$ by convexity.
	Proposition \ref{sup norm est} gives us the sup norm estimate:
	\begin{align}\label{S3 Linfty convergence}
		\|v(t,\,\cdot\,)-M\,Q_{t}^{\sigma}\|_{L^{\infty}(\mathbb{X})}\,
		\le\,
		\|f\|_{L^{1}(\mathbb{X})}\,\|Q_{t}^{\sigma}\|_{L^{\infty}(\mathbb{X})}\,
		+\,
		|M|\,\|Q_{t}^{\sigma}\|_{L^{\infty}(\mathbb{X})}\,
		\lesssim\,t^{\sigma-\frac{\ell}{2}-\frac{1}{2}-|\Sigma_r^+|}e^{-|\rho|t}
	\end{align}
	for $t$ large and for all $f\in{L^{1}(\mathbb{X})}$.
	Notice that such an estimate holds without the bi-$K$-invariance assumption.
	By convexity, we obtain the following estimates in the $L^{p}(\mathbb{X})$
	setting.
	\begin{corollary}
		Under the assumptions of \cref{S1 Main thm 1}, we have
		\begin{align}\label{S3 Lp convergence}
			\|v(t,\,\cdot\,)-M\,Q_{t}^{\sigma}\|_{L^{p}(\mathbb{X})}\,
			=\,\mathrm{o}\big(t^{-\frac{1}{p'}(\sigma-\frac{\ell}{2}-\frac{1}{2}-|\Sigma_r^+|)}
			e^{-\frac{|\rho|t}{p'}}\big)
			\qquad\textnormal{as}\quad\,t\rightarrow+\infty
		\end{align}
		for all $1<p<\infty$.
	\end{corollary}
	
	\begin{remark}
		The sup norm estimate \eqref{S3 Linfty convergence} is weaker 
		compared to the results in the Euclidean setting. More precisely, on $\mathbb{R}^n$, the Poisson semigroup ($\sigma=1/2$) satisfies the strong convergence \eqref{S1 Linf R Poisson} (recall that $\|\mathcal{Q}_t\|_{L^{\infty}(\mathbb{R}^n)}\asymp t^{-n}$). However, this is not true on noncompact symmetric spaces.
		Indeed, in the lines of \cite[Remark 3.6]{APZ2023}, consider the Poisson kernel $Q_t^{1/2}$ as well as a ``delayed'' Poisson kernel $Q_{t+t'}^{1/2}$ for some $t'>0$ to be determined 
		later. Recall that $\nu=\ell+2|\Sigma_r^+|$. Then 
		\begin{align*}
			t^{\frac{\nu}{2}}e^{|\rho|t}\,
			\|Q_{t+t'}^{1/2}-Q_{t}^{1/2}\|_{L^{\infty}(\mathbb{X})}\,
			\ge\,t^{\frac{\nu}{2}}e^{|\rho|t}\,\big(Q_{t}^{1/2}(eK)-Q_{t+t'}^{1/2}(eK)\big)
		\end{align*}
		since $Q_{t}^{1/2}(eK)$ is decreasing in $t$, as seen by the subordination formula.
		According to \eqref{P-bound}, there exists a constant $C\ge1$ such that
		\begin{align*}
			t^{\frac{\nu}{2}}\,e^{|\rho|t}\,\big(Q_{t}^{1/2}(eK)-Q_{t+t'}^{1/2}(eK)\big)
			&\ge\,
			t^{\frac{\nu}{2}}\,e^{|\rho|t}\,
			\big\lbrace{
				\tfrac{1}{C}\,t^{-\frac{\nu}{2}}\,e^{-|\rho|t}\,
				-\,
				C(t+t')^{-\frac{\nu}{2}}\,e^{-|\rho|(t+t')}
			}\big\rbrace\\[5pt]
			&=\,
			C^{-1}\,-\,C\big(\tfrac{t}{t+t'}\big)^{\frac{\nu}{2}}\,e^{-|\rho|t'}
			\ge\,\tfrac{1}{2C}\,,
		\end{align*}
		provided that $t'>\frac{2\ln{C}+\ln2}{|\rho|}$.
		Hence
		\begin{align*}
			t^{\frac{\nu}{2}}\,e^{|\rho|t}\,
			\|Q_{t+t'}^{1/2}-Q_{t}^{1/2}\|_{L^{\infty}(\mathbb{X})}\,
			\centernot\longrightarrow\,0
			\qquad\textnormal{as}\quad\,t\rightarrow +\infty.
		\end{align*}
	\end{remark}

	\section{Asymptotic convergence associated with 
		the extension problem for the distinguished Laplacian}\label{Section.5 Distinguished}
	
	Let $S=N(\exp{\mathfrak{a}})=(\exp{\mathfrak{a}})N$ be the solvable group
	occurring in the Iwasawa decomposition $G=N(\exp{\mathfrak{a}})K$. 
	Then $S$ is identifiable, as a manifold, with the symmetric space 
	$\mathbb{X}=G/K$. The distinguished Laplacian $\widetilde{\Delta}$ on $S$ is 
	given by the conjugation of the shifted Laplace-Beltrami operator
	$\Delta+|\rho|^{2}$ on $\mathbb{X}$:
	\begin{align}\label{S4 dist laplacian}
		\widetilde{\Delta}\,
		=\,\widetilde{\delta}^{\frac12}\circ(\Delta+|\rho|^{2})
		\circ\widetilde{\delta}^{-\frac12}
	\end{align}
	where the modular function $\widetilde{\delta}$ of $S$ is defined by
	\begin{align*}
		\widetilde{\delta}(g)\,
		=\,\widetilde{\delta}(n(\exp{A}))\,
		=\,e^{-2\langle{\rho,A}\rangle}
		\qquad\forall\,g\in{S}.
	\end{align*}
	Here $n=n(g)$ and $A=A(g)$ denotes respectively the $N$-component and the
	$\mathfrak{a}$-component of $g$ in the Iwasawa decomposition.
	The distinguished Laplacian $\widetilde{\Delta}$ is left-$S$-invariant and
	self-adjoint with respect to the right-invariant Haar measure on $S$:
	\begin{align*}
		\int_{S}\textrm{d}_{r}{g}\,f(g)
		=\,\int_{N}\diff{n}\int_{\mathfrak{a}}\diff{A}\,f(n(\exp{A}))\,
		=\,\int_{\mathfrak{a}}\diff{A}\,e^{2\langle{\rho,A}\rangle}
		\int_{N}\diff{n}f((\exp{A})n).
	\end{align*}
	The connection between the measures on $S$
	and the unimodular Haar measure on $G$ is given as follows,
	\begin{align}\label{S4 drdldg}
		\int_{S}\textrm{d}_{r}{g}\,f(g)\,
		=\,\int_{G}\diff{g}\,e^{2\langle{\rho,A(g)}\rangle}f(g)
		\qquad\textnormal{and}\qquad
		\int_{S}\textrm{d}_{\ell}{g}\,f(g)\,
		=\,\int_{G}\diff{g}\,f(g).
	\end{align}
	
	Recall the heat equation associated with the distinguished Laplacian:
	\begin{align}\label{HE S}
		\partial_{t}\widetilde{v}(t,g)\,
		=\,\widetilde{\Delta}_{g}\widetilde{v}(t,g),
		\qquad
		\widetilde{v}(0,g)\,=\,\widetilde{v}_{0}(g),
	\end{align}
	where the corresponding heat kernel is given by 
	$\widetilde{h}_{t}=\widetilde{\delta}^{1/2}\,e^{|\rho|^{2}t}h_{t}$
	in the sense that
	\begin{align*}
		(e^{t\widetilde{\Delta}\,}f)(g)\,
		=\,(f*\widetilde{h}_{t})(g)\,
		=\,\int_{S}\textrm{d}_{\ell}{y}\,f(y)\,\widetilde{h}_{t}(y^{-1}g)\,
		=\,\int_{S}\textrm{d}_{r}{y}\,f(gy^{-1})\,\widetilde{h}_{t}(y).
	\end{align*}
	Here, we still denote by $*$ the convolution product on $S$ or on $G$. We refer 
	to \cite{Bou1983,CGGM1991} for more details about the distinguished Laplacian.

	The extension problem for $\widetilde{\Delta}$ on $S$ writes
	\begin{align}\label{extensionS}
		\begin{cases}
			\widetilde{\Delta} \widetilde{v} - \frac{(1-2\sigma)}{t} \partial_{t} \widetilde{v} - \partial^{2}_{tt} \widetilde{v}=0, \qquad  \; t>0, \\
			\widetilde{u}(\, \cdot \, ,0)=\widetilde{v}_{0}
		\end{cases} 
	\end{align}
	which admits a fundamental solution $\widetilde{Q}_t^{\sigma}$. The use of the latter is justified by the arguments in \cite{ST2010}  discussed in Section \ref{Section.3 CS}, and more precisely by \eqref{L2heat}, by the fact that
	\[
	\partial_t\widetilde{h}_t=\widetilde{\delta}^{1/2}\,\partial_t (e^{(|\rho|^2t}h_t)=\widetilde{\delta}^{1/2}\,e^{|\rho|^2t}\big(|\rho|^2\,h_t+\partial_t h_t\big), \quad t\in (0,1),
	\] 
	and by the observation that if $f\in L^2(\mathbb{X})$ and $\widetilde{f}=\widetilde{\delta}^{1/2}f$, then by~\eqref{S4 drdldg} we have
	\[
	\int_{S}\diff_r g\,\widetilde{f}^2(g)  =\int_{S}\diff_r g\, e^{-2\langle \rho, A(g)\rangle} f^2(g)=\int_{G} \diff{g}\,f^2(g) .
	\]
	Therefore, by the Stinga-Torrea subordination formula we have
	\begin{align}
		\widetilde{Q}_t^{\sigma}(g)&=\frac{t^{2\sigma}}{2^{2\sigma}\Gamma(\sigma)}\int_{0}^{+\infty} \frac{\diff u}{u^{1+\sigma}}\,\widetilde{h}_{u}(g)\,e^{-\frac{t^2}{4u}} \label{Q kernel S0} \\
		&=\widetilde{\delta}^{\frac{1}{2}}(g)\frac{t^{2\sigma}}{2^{2\sigma}\Gamma(\sigma)}\int_{0}^{+\infty}\frac{\diff u}{u^{1+\sigma}} \,e^{|\rho|^2u}\,h_u(g)\,e^{-\frac{t^2}{4u}}=:\widetilde{\delta}^{\frac{1}{2}}(g) Q^{\sigma, 0}(g). \label{Q kernel S}
	\end{align}

	\begin{remark}
		Notice that $\widetilde{Q}_{t}^{\sigma}(g)\textrm{d}_{r}{g}$ is a probability measure
		on $S$. Indeed, this follows from the subordination formula \eqref{Q kernel S0}  and the fact that $\int_{S}\textrm{d}_{r}{g}\,\widetilde{h}_t(g)=1$.
	\end{remark}
	
	The first subsection is devoted to determine the critical region 
	where the  kernel $\widetilde{Q}_{t}^{\sigma}$ concentrates. 
	In the next two subsections, we study respectively the $L^{1}$ 
	and the $L^{\infty}$ asymptotic convergences of solutions to \eqref{HE S}
	with compactly supported initial data (no bi-$K$-invariance required). 
	We discuss the same questions for other initial data in the last subsection.
	
	\subsection{Asymptotic concentration of the fractional Poisson kernel associated to the distinguished Laplacian}
	We first give large time asymptotics for the distinguished extension kernel. More precisely, we prove the following upper and lower bounds.
	
	\begin{proposition}\label{kernel bounds S}
		The  fractional Poisson kernel $\widetilde{Q}_t^\sigma$ on $S$, $0<\sigma<1$, associated with the distinguished Laplacian, satisfies the upper and lower bounds 
		\begin{equation}\label{P-bound dist}
			\widetilde{Q}_t^{\sigma} (g)\asymp  \widetilde{\delta}^{\frac{1}{2}}(g)\, \varphi_0(g)\,
			t^{2\sigma}\, (t+|g^{+}|)^{-\ell-2|\Sigma_r^+|-2\sigma},
		\end{equation}
		if $t^2+|g^{+}|^2 \geq 1$.
	\end{proposition}

	\begin{proof}
		Recall that by the subordination formula \eqref{Q kernel S}, we may write
		$$	\frac{2^{2\sigma}\Gamma(\sigma)}{t^{2\sigma}}\widetilde{\delta}^{-\frac{1}{2}}(g)\,\widetilde{Q}_t^{\sigma} (g)=\int_{0}^{t^2+|g|^2} \frac{\diff u}{u^{1+\sigma}}\,h_u(g)\,e^{|\rho|^2 u}\,e^{-\frac{t^2}{4u}}
		+\int_{t^2+|g|^2}^{+\infty} \frac{\diff u}{u^{1+\sigma}}\,h_u(g)\,e^{|\rho|^2 u}\,e^{-\frac{t^2}{4u}},$$
		and observe that due to \eqref{S2 heat kernel} we have
		$$u^{-1-\sigma} \,h_u(g)\,e^{|\rho|^2 u}\,e^{-\frac{t^2}{4u}}\asymp\varphi_{0}(g)\,u^{-\frac{n}{2}-1-\sigma}\,
		\Big\lbrace{
			\prod_{\alpha\in\Sigma_{r}^{+}}
			(1+u+\langle{\alpha,g^{+}}\rangle)^{\frac{m_{\alpha}+m_{2\alpha}}{2}-1}
		}\Big\rbrace\,
		e^{-\frac{t^2+|g^{+}|^{2}}{4u}}, \quad u>0.$$
		
		Let $t^2+|g^{+}|^2 \geq 1$. As far as upper bounds are concerned,  the claim follows in the first interval $(0, t^2+|g^{+}|^2]$ by estimating $1+u+\langle \alpha, g^{+}\rangle\lesssim t^2+|g^{+}|^2$. For the second interval  $(t^2+|g^{+}|^2, +\infty)$ we estimate $1+u+\langle \alpha, g^{+}\rangle\lesssim u$, and take into account that $\sum_{\alpha\in\Sigma_{r}^{+}} (m_{\alpha}+m_{2\alpha})=n-\ell.$
		
		The lower bound follows writing 
		\[
		1+u+\langle \alpha, g^+\rangle \geq  u
		\]
		and integrating over $(t^2+|g|^2, +\infty)$. We omit the details.
	\end{proof}	

	Recall that the extension kernel $Q_{t}^{\sigma}$ associated with the Laplace-Beltrami operator
	concentrates in $K(\exp\Omega_{t})K$, where $\Omega_{t}$ is described in Proposition \ref{prop critical region}.
	The following proposition shows that the kernel $\widetilde{Q}_{t}^{\sigma}$ 
	(associated with the extension problem for the distinguished Laplacian) concentrates in a different region.
	
	\begin{proposition}\label{prop critical region S} Let $0<\varepsilon<1$. Consider in $\mathfrak{a}$ the annulus
		\begin{equation} 
			B(0,t^{1+\varepsilon})\smallsetminus B(0,t^{1-\varepsilon}).
		\end{equation}
		Then, the fractional Poisson kernel associated with the distinguished Laplacian on $S$
		concentrates asymptotically in $K(\exp\widetilde{\Omega}_{t})K$.	In other words,
		\begin{align*}
			\lim_{t\rightarrow +\infty}
			\int_{g\in{S}\,\textrm{s.t.}\,
				g^{+}\in\overline{\mathfrak{a}^{+}}
				\smallsetminus\widetilde{\Omega}_{t}}
			\textrm{d}_{r}{g}\,\widetilde{Q}_{t}^{\sigma}(g)\,
			=\,0
		\end{align*}
		where $g^{+}$ denotes the middle component of $g$ in the Cartan decomposition.
	\end{proposition}
	\begin{proof}
		By using \eqref{S4 drdldg} and \eqref{Q kernel S},  write
		\begin{align*}
			I(t)={\int_{g\in{S}\,\textrm{ s.t.}\,
					g^{+}\in\overline{\mathfrak{a}^{+}}
					\smallsetminus\widetilde{\Omega}_{t}}} 
			\textrm{d}_{r}{g}\,\widetilde{Q}_{t}^{\sigma}(g)
			=\int_{G\smallsetminus K(\exp\widetilde{\Omega}_{t})K}\diff{g}\,
			e^{\langle{\rho,A(g)}\rangle}Q_t^{\sigma,0}(g).
		\end{align*}

		Since $Q_t^{\sigma,0}$ is bi-$K$-invariant on $G$, writing
		$\diff{k}$ for the normalized Haar measure on the compact group $K$, and using Proposition \ref{kernel bounds S} we have
		\begin{align*}
			I(t)\,
			&\asymp t^{2\sigma}\,\int_{G\smallsetminus K(\exp\widetilde{\Omega}_{t})K}\diff{g}\,\varphi_0(g)\,
			(t+|g|)^{-\ell-2|\Sigma_r^+|-2\sigma}
			\int_{K}\diff{k}\,e^{\langle{\rho,A(kg)}\rangle}\\
			&\asymp t^{2\sigma}\,
			\int_{G\smallsetminus K(\exp\widetilde{\Omega}_{t})K}\diff{g}\,\varphi_0(g)^2\,
			(t+|g|)^{-\ell-2|\Sigma_r^+|-2\sigma}.
		\end{align*}
		According to the Cartan decomposition, and to the estimates
		\eqref{S2 estimate of delta} and \eqref{S2 global estimate phi0}, we obtain
		\begin{align}\label{S4 global estimate I}
			I(t)\,
			&\asymp \, t^{2\sigma}\,
			\int_{\widetilde{\Omega}_{t}^{c}}\diff{g^{+}}\,
			\delta(g^{+})\,(t+|g^{+}|)^{-\ell-2|\Sigma_r^+|-2\sigma}\,\varphi_{0}(\exp{g^{+}})^2
			\notag\\[5pt]
			&\lesssim\,t^{2\sigma}\,
			\int_{\widetilde{\Omega}_{t}^{c}}\diff{g^{+}}\,
			(t+|g^{+}|)^{-\ell-2|\Sigma_r^+|-2\sigma}\,
			(1+|g^{+}|)^{2|\Sigma_r^{+}|}.
		\end{align}

		Next, let us study the right-hand side of \eqref{S4 global estimate I} 
		outside $\widetilde{\Omega}_{t}$. On the one hand, if $|g^{+}|<t^{1-\varepsilon}$ then $t+|g^{+}|\asymp t$, so
		\begin{align*}
			t^{2\sigma}	\int_{|g^{+}|<t^{1-\varepsilon}}\diff{g^{+}}\,
			(t+|g^{+}|)^{-\ell-2|\Sigma_r^+|-2\sigma}
			(1+|g^{+}|)^{2|\Sigma_r^{+}|} &\lesssim	t^{2\sigma} \int_{0}^{t^{1-\varepsilon}}\diff r\, r^{\ell-1} t^{-\ell-2|\Sigma_r^+|-2\sigma}\,(1+r)^{2|\Sigma_r^{+}|}\\
			&\lesssim t^{-\varepsilon(\ell+2|\Sigma_r^{+}|)}.
		\end{align*}

		On the other hand, if $|g^{+}|>t^{1+\varepsilon}$ then $t+|g^{+}|\asymp |g^{+}|$, so we have
		\begin{align*}
			t^{2\sigma}	\int_{|g^{+}|>t^{1+\varepsilon}}\diff{g^{+}}\,
			(t+|g^{+}|)^{-\ell-2|\Sigma_r^+|-2\sigma}
			(1+|g^{+}|)^{2|\Sigma_r^{+}|}
			&\lesssim	t^{2\sigma} \int_{t^{1+\varepsilon}}^{+\infty}\diff r\, r^{\ell-1} r^{-\ell-2|\Sigma_r^+|-2\sigma}\,(1+r)^{2|\Sigma_r^{+}|}\\
			&\lesssim t^{-\varepsilon\sigma}.
		\end{align*}
		In other words, we have proved that $I(t)=\textrm{O}(t^{-\sigma \varepsilon})$, therefore the Poisson  kernel $\widetilde{Q}_{t}^{\sigma}$ associated with the 
		distinguished Laplacian on $S$ concentrates asymptotically in 
		$K(\exp\widetilde{\Omega}_{t})K$.
	\end{proof}
	
	\begin{remark} The critical region for the fractional Poisson kernel associated with the distinguished Laplacian is similar to that of its Euclidean counterpart.
	\end{remark}	
	
	We now obtain  precise long-time asymptotics of the kernel 
	$$Q_t^{\sigma,0}(g)=\frac{t^{2\sigma}}{2^{2\sigma}\Gamma(\sigma)}\int_{0}^{+\infty}\frac{\diff u}{u^{1+\sigma}} \,e^{|\rho|^2u}\,h_u(g)\,e^{-\frac{t^2}{4u}},$$ 
	which are crucial for our proof, by a slightly more general result.
	
	\begin{theorem} \label{asymp kernel S}
		Let $\sigma \in (0, 1)$ and  $g\in S$ such that $g^{+}\in \widetilde{\Omega_{t}}$. Then, as $t+|g^{+}|\rightarrow+\infty$, we have
		\begin{align}
			Q_t^{\sigma,0}(g) \sim  \widetilde{C}(\sigma)\,
			t^{2\sigma}\, \varphi_0(\exp g^{+}) \left(t^2+|g^{+}|^2\right)^{-\frac{\ell}{2}-|\Sigma_r^+|-\sigma},
		\end{align}
		where the constant is 
		$$\widetilde{C}(\sigma)=\frac{1}{ \Gamma(\sigma)}C_0 \, 2^{\ell+|\Sigma_r^{+}|} \pi^{\frac{\ell}{2}}\Gamma \left( \frac{\ell}{2}+|\Sigma_r^{+}|+\sigma\right) \bm{\pi}(\widetilde \rho)\bm{b}(0)^{-2},$$
		with $\widetilde{\rho}=\frac{1}{2}\sum_{\alpha \in \Sigma_r^{+}}\alpha$ and $C_0=2^{n-\ell}/(2\pi)^{\ell}|K/\mathbb{M}|$.
	\end{theorem}	
	
	\begin{proof}
		The proof follows arguments for the asymptotics of the Poisson kernel ($\sigma=1/2$) in \cite[Theorem 5.3.1]{AnJi1999}. 
		
		Let $0<\eta<1$ and $t^2+|g^{+}|^2>1$. In view of the subordination formula \eqref{Q kernel S}, let us split 
		\begin{align*}
			Q_t^{\sigma,0}(g) &= \frac{t^{2\sigma}}{2^{2\sigma}\Gamma(\sigma)}\int_{0}^{+\infty}\frac{\diff u}{u^{1+\sigma}} \,e^{|\rho|^2u}\,h_u(g)\,e^{-\frac{t^2}{4u}}\\
			&=\frac{t^{2\sigma}}{ 2^{2\sigma} \Gamma(\sigma)}  \{J_1+J_2+J_3\},
		\end{align*}
		where the quantities $J_1, J_2$ and $J_3$ are defined by the integration  over the intervals $\big[0, (t^2+|g|^2)^{1-\eta} \big), \big[(t^2+|g|^2)^{1-\eta}, (t^2+|g|^2)^{1+\eta}\big)$ and $\big[(t^2+|g|^2)^{1+\eta}, \infty\big)$ respectively. 
		
		We claim that the main contribution comes from the middle integral $J_2$. Indeed,  for the first integral $J_1$, working as in Proposition \ref{kernel bounds S}, we get 	$J_1 =\textrm{O}((t^2+|g|^2)^{-\infty}\varphi_0(g))$ while, similarly, for the third integral we get $J_3=\textrm{O} ( (t^2+|g|^2)^{-(1+\eta)(\frac{\ell}{2}+|\Sigma_r^+|+\sigma)}\varphi_{0}(g))$.

		We now consider $J_2$. Define
		\begin{align}\label{h' kernel}
			h'(t,g^{+})= t^{\frac{\ell}{2}+|\Sigma_{r}^{+}|}\,
			\varphi_{0}(\exp{g^{+}})^{-1}\,e^{|\rho|^{2}t+\frac{|g^{+}|^{2}}{4t}}\, h_t(\exp{g^{+}}), \quad t>0, \, g^{+}\in \overline{\mathfrak{a}^{+}}.
		\end{align}
		Then,  by \eqref{h' kernel} we have
		\begin{align*}
			J_2&= \int_{ (t^2+|g|^2)^{1-\eta}}^{ (t^2+|g|^2)^{1+\eta}} \frac{\diff u}{u^{1+\sigma}}\, e^{|\rho|^2u}h_u(\exp g^{+})\,e^{-\frac{t^2}{4u}}\\
			&=\varphi_0(\exp g^{+})\int_{ (t^2+|g|^2)^{1-\eta}}^{ (t^2+|g|^2)^{1+\eta}}\diff u \, u^{-\frac{\ell}{2}-|\Sigma_r^{+}|-\sigma-1}e^{-\frac{t^2+|g^{+}|^2}{4u}}h'(u, g^{+})\\
			&=2^{\ell+2|\Sigma_r^{+}|+2\sigma}(t^2+|g^{+}|^2)^{-\frac{\ell}{2}-|\Sigma_r^{+}|-\sigma}\int_{ \frac{1}{4}(t^2+|g|^2)^{-\eta}}^{ \frac{1}{4}(t^2+|g|^2)^{\eta}}\diff u \, u^{\frac{\ell}{2}+|\Sigma_r^{+}|+\sigma-1}e^{-u}h'\left(\frac{t^2+|g^{+}|^2}{4u}, g^{+}\right).
		\end{align*}
		Since
		$$h'\left(\frac{t^2+|g^{+}|^2}{4u}, g^{+}\right) \longrightarrow C_2\,\bm{b}(0)^{-1},$$
		uniformly 	as $t^2+|g^{+}|^2\rightarrow +\infty$, by \eqref{AJ heat asymp} and contradiction (see \cite[p.1086]{AnJi1999}), the Laplace method we obtain that the last integral tends to 
		$$C_2\,\bm{b}(0)^{-1}\Gamma \left( \frac{\ell}{2}+|\Sigma_r^{+}|+\sigma\right).$$
		
		Since $J_1, J_3$ are very small compared to $J_2$ for $t$ large, substituting the value of $C_2$ (see \eqref{AJ heat asymp}) we finally get the claimed asymptotics.
	\end{proof}
	
	
	\subsection{Heat asymptotics in $L^1$ for compactly supported initial data}\label{Subsection L1 S}
	In this subsection, we investigate the long-time asymptotic convergence in 
	$L^1(S)$ of solutions to the Cauchy problem \eqref{HE S}, where the initial data 
	$\widetilde{v}_{0}$ is assumed continuous and compactly supported in $B(eK,\xi)$. 
	Let
	$\widetilde{\varphi}_{0}=\widetilde{\delta}^{\frac12}\varphi_{0}$
	be the modified ground spherical function. The mass function is defined by
	\begin{align}\label{S4 mass G}
		\widetilde{M}(g)\,
		=\,\frac{(\widetilde{v}_{0}*\widetilde{\varphi}_{0})(g)
		}{\widetilde{\varphi}_{0}(g)}
		\qquad\forall\,g\in{S}.
	\end{align}
	By using the fact that the modular function $\widetilde{\delta}$ is a character
	on $S$, we can also write the mass as
	\begin{align}\label{S4 mass S}
		\widetilde{M}(g)\,
		=\,\tfrac{1}{\widetilde{\delta}(g)^{\frac12}\,\varphi_{0}(g)}\,
		\int_{S}\textrm{d}_{\ell}{y}\,
		v_{0}(gK)\,
		\underbrace{\vphantom{\Big|}
			\widetilde{\delta}(y)^{\frac12}
			\widetilde{\delta}(y^{-1}g)^{\frac12}
		}_{\widetilde{\delta}(g)^{\frac12}}\,
		\varphi_{0}(y^{-1}g)
		=\,\frac{(v_{0}*\varphi_{0})(g)}{\varphi_{0}(g)}
	\end{align}
	where $v_{0}(gK)=\widetilde{\delta}(g)^{-\frac{1}{2}}\widetilde{v}_{0}(g)$
	is a right $K$-invariant function on $G$, 
	with compact support $(\supp\widetilde{v}_{0})K$.

	The following properties of the mass function were already observed in \cite[Remarks 4.5 and 4.6]{APZ2023}:
	
	\begin{remark}
		1. If $\widetilde{v}_{0}\in\mathcal{C}_{c}(S)$,
		then the mass function $\widetilde{M}$ is bounded. This follows from the fact that 
		\begin{equation}\label{Harnack}
			\frac{\varphi_{0}(y^{-1}g)}{\varphi_0(g)}\leq C(\xi) \quad \text{if} \quad |y|<\xi.
		\end{equation}
		
		2.	The mass function $\widetilde{M}$ is a constant 
		\noindent		if $v_{0}$ is bi-$K$-invariant and 
		$\widetilde{v}_{0}=\widetilde{\delta}^{\frac12}v_{0}$ belongs to $L^{1}(S)$:
		\begin{align*}
			\widetilde{M}\,
			=\,\int_{G}\diff{y}\,v_{0}(y)\,\varphi_{0}(y)
			=\,\mathcal{H}v_{0}(0).
		\end{align*}
	\end{remark}
	
	The following lemma plays a key role in the proof of Theorem \ref{S1 Main thm 2}.

	\begin{lemma}\label{S4 Lemma ratios difference}
		For  bounded $y\in{G}$ and for all $g$ in the critical region
		$K(\exp\widetilde{\Omega}_{t})K$, the following asymptotic behavior holds:
		\begin{align*}
			\frac{Q_{t}^{\sigma, 0}(y^{-1}g)}{Q_{t}^{\sigma, 0}(g)}
			-\frac{\varphi_{0}(y^{-1}g)}{\varphi_{0}(g)}\,
			=\,\mathrm{O}\Big(t^{-1+\varepsilon(\nu+2\sigma-1)}\Big)
			\qquad\textnormal{as}\,\,\,t\rightarrow +\infty,
		\end{align*}
		for $0<\varepsilon<1/(\nu+2\sigma-1)$. Here, $\nu=\ell+2|\Sigma_r^{+}|.$
	\end{lemma}
	
	\begin{proof}
		Assume that $|y|\le\xi$ for some positive constant $\xi$.
		Recall that for every $H\in\widetilde{\Omega}_{t}$, we have 
		$t^{1-\varepsilon}\le|H|\le t^{1+\varepsilon}$.
		Notice also that 
		\begin{align*}
			|(y^{-1}g)^{+}-g^{+}|\,\le\,|y|\,<\,\xi
		\end{align*}
		according to \eqref{dist flat}. Then, for $t$ large enough, we deduce the following estimates:
		\begin{align*}
			\begin{cases}
				|(y^{-1}g)^{+}|\,\le\,|g^{+}|+\xi\,
				<\,t^{1+\varepsilon}+\xi\,<2\,t^{1+\varepsilon},\\[5pt]
				|(y^{-1}g)^{+}|\,\ge\,|g^{+}|-\xi\,
				>\,t^{1-\varepsilon}-\xi\,>\frac{1}{2}\,t^{1-\varepsilon},\\[5pt]
			\end{cases}
		\end{align*}
		In other words, we obtain
		\begin{align*}
			y^{-1}g\,\in\,K(\exp{\widetilde{\Omega}_{t}'})K
			\qquad\forall\,g\in{K(\exp\widetilde{\Omega}_{t})K},\,\,
			\forall\,|y|<\xi,
		\end{align*}
		where
		\begin{align*}
			\widetilde{\Omega}_{t}'\,
			=\,\big\lbrace{g\in S:\; 
				\frac{1}{2}\,t^{1-\varepsilon}\le|g^{+}|\le 2\,t^{1+\varepsilon}
			}\big\rbrace.
		\end{align*}

		Thus the asymptotics of Theorem \ref{asymp kernel S} yield
		\begin{align*}
			\frac{\widetilde{Q}_{t}^{\sigma}(y^{-1}g)}{\widetilde{Q}_{t}^{\sigma}(g)}	-\frac{\varphi_{0}(y^{-1}g)}{\varphi_{0}(g)}\,
			\sim\, \frac{\varphi_{0}(y^{-1}g)}{\varphi_{0}(g)}\left(   \frac{\left(t^2+|g^{+}|^2\right)^{\frac{\ell}{2}+|\Sigma_r^+|+\sigma}}{\left(t^2+|(y^{-1}g)^{+}|^2\right)^{\frac{\ell}{2}+|\Sigma_r^+|+\sigma}} -1\right) ,
		\end{align*}
		On the one hand, the quotient of the ground spherical functions is bounded by the local Harnack inequality \eqref{Harnack}. On the other hand, using \eqref{MVTk} for $r=|g^{+}|$, $s=|(y^{-1}g)^{+}|$ and $k=\frac{\ell}{2}+|\Sigma_r^{+}|+\sigma=\frac{\nu}{2}+\sigma > 1$ (in the notation of \eqref{MVTk}, we now have $r_0\lesssim t^{1+\varepsilon}$) and the trivial inequality $t^2+s^2\geq t^2$, we get altogether
		\begin{align*}
			\frac{Q_{t}^{\sigma,0}(y^{-1}g)}{Q_{t}^{\sigma,0}(g)}
			-\frac{\varphi_{0}(y^{-1}g)}{\varphi_{0}(g)}\,
			=\,\mathrm{O}\Big(t^{-1+\varepsilon(\nu+2\sigma-1)}\Big)
			\qquad\forall\,g\in{K(\exp\widetilde{\Omega}_{t})K},\,\,
			\forall\,|y|<\xi.
		\end{align*}
	\end{proof}

	\vspace{5pt}
	Now, let us prove the first part of \cref{S1 Main thm 2}. The arguments follow those of \cite{APZ2023} once Lemma \ref{S4 Lemma ratios difference} is at hand, but we include them for the reader's convenience.
	\begin{proof}[Proof of \eqref{S1 L1 disting} in \cref{S1 Main thm 2}.]
		
		By using
		\begin{align*}
			(\widetilde{v}_{0}*\widetilde{\varphi}_{0})(g)\,
			=\,\int_{S}\textrm{d}_{\ell}{y}\,
			v_{0}(yK)
			\underbrace{\vphantom{\Big|}
				\widetilde{\delta}(y)^{\frac12}\,\widetilde{\delta}(y^{-1}g)^{\frac12}
			}_{\widetilde{\delta}(g)^{\frac12}}
			\varphi_{0}(y^{-1}g)
			=\,\widetilde{\delta}(g)^{\frac{1}{2}}
			(v_{0}*\varphi_{0})(gK),
		\end{align*}
		and the fact that $\widetilde{Q}_t^{\sigma}=\widetilde{\delta}^{\frac{1}{2}}Q_t^{\sigma,0}$, let us write the solution $\widetilde{v}$ to \eqref{extensionS} as
		\begin{align*}
			\widetilde{v}(t,g)\,
			=\,(\widetilde{v}_{0}*\widetilde{Q}_{t}^{\sigma})(g)\,
			=\,\widetilde{\delta}(g)^{\frac{1}{2}}\,
			(v_{0}*Q_{t}^{\sigma, 0})(g).
		\end{align*}

		We aim to study the difference
		\begin{align}\label{S4 difference}
			\widetilde{v}(t,g)-\widetilde{M}(g)\widetilde{Q}_{t}^{\sigma}(g)\,
			&=\,\widetilde{Q}_{t}^{\sigma}(g)\,
			\frac{(v_{0}*Q_{t}^{\sigma,0})(g)}{Q_{t}^{\sigma,0}(g)}
			\,-\,
			\widetilde{Q}_{t}^{\sigma}(g)\,
			\frac{(v_{0}*\varphi_{0})(g)}{\varphi_{0}(g)}
			\notag\\[5pt]
			&=\,\widetilde{Q}_{t}^{\sigma}(g)
			\int_{G}\diff{y}\,v_{0}(yK)\,
			\Big\lbrace{
				\frac{Q_{t}^{\sigma,0}(y^{-1}g)}{Q_{t}^{\sigma,0}(g)}
				-\frac{\varphi_{0}(y^{-1}g)}{\varphi_{0}(g)}
			}\Big\rbrace.
		\end{align}
		According to the previous lemma, we have
		\begin{align*}
			\frac{Q_{t}^{\sigma,0}(y^{-1}g)}{Q_{t}^{\sigma,0}(g)}
			-\frac{\varphi_{0}(y^{-1}g)}{\varphi_{0}(g)}\,
			=\,\mathrm{O}\Big(t^{-1+\varepsilon(\nu+2\sigma-1)}\Big)
			\qquad\forall\,g\in{K(\exp\widetilde{\Omega}_{t})K},\,\,
			\forall\,y\in\supp{v_{0},}
		\end{align*}
		and therefore the integral of
		$\widetilde{v}(t,\cdot)-\widetilde{M}\,\widetilde{Q}_{t}^{\sigma}$ 
		over the critical region
		\begin{align*}
			\int_{S\cap{K(\exp\widetilde{\Omega}_{t})K}}
			\textrm{d}_{r}{g}\,    
			|\widetilde{v}(t,g)-\widetilde{M}(g)\widetilde{Q}_{t}^{\sigma}(g)|\,
			\lesssim\,
			t^{-1+\varepsilon(\nu+2\sigma-1)}\,
			\underbrace{\vphantom{\Big|}
				\int_{S}\textrm{d}_{r}{g}\,\widetilde{Q}_{t}^{\sigma}(g)
			}_{1}\,
			\underbrace{\vphantom{\Big|}
				\int_{G}\diff{y}\,|v_{0}(yK)|
			}_{\const}
		\end{align*}
		tends asymptotically to $0$. Finally, we claim that the integral
		\begin{align*}
			\int_{S\smallsetminus{K(\exp\widetilde{\Omega}_{t})K}}\textrm{d}_{r}{g}\,
			|\widetilde{v}(t,g)-\widetilde{M}(g)\widetilde{Q}_{t}^{\sigma}(g)|\,
			&\le\,
			\int_{S\smallsetminus{K(\exp\widetilde{\Omega}_{t})K}}\textrm{d}_{r}{g}\,
			|\widetilde{v}(t,g)|\\[5pt]
			&+\,
			\int_{S\smallsetminus{K(\exp\widetilde{\Omega}_{t})K}}\textrm{d}_{r}{g}\,    
			|\widetilde{M}(g)|\widetilde{Q}_{t}^{\sigma}(g)
		\end{align*}
		tends also to $0$. On the one hand, we know that $\widetilde{M}$ is bounded
		and that the kernel $\widetilde{Q}_{t}^{\sigma}$ asymptotically concentrates in 
		$K(\exp\widetilde{\Omega}_{t})K$, hence
		\begin{align*}
			\int_{S\smallsetminus{K(\exp\widetilde{\Omega}_{t})K}}\textrm{d}_{r}{g}\,
			|\widetilde{M}(g)|\widetilde{Q}_{t}^{\sigma}(g)\,
			\longrightarrow\,0
		\end{align*}
		as $t\rightarrow +\infty$. On the other hand, notice that for all
		$y\in\supp{v_{0}}$ and for all $g\in{G}$ such that
		$g^{+}\notin\widetilde{\Omega}_{t}$, using the triangle inequality one can show that
		\begin{align}\label{S4 Omega''}
			(y^{-1}g)^{+}\,\notin\,\widetilde{\Omega}_{t}''
			=\,
			\big\lbrace{
				H\in\overline{\mathfrak{a}^{+}}\,|\,
				2\,t^{1-\varepsilon}\le|H|\le\tfrac{1}{2}\,t^{1+\varepsilon}\big\rbrace}.
		\end{align}

		Hence
		\begin{align*}
			\int_{S\smallsetminus{K(\exp\widetilde{\Omega}_{t})K}}\textrm{d}_{r}{g}\,
			|\widetilde{v}(t,g)|\,
			&\le\,\int_{G}\diff{y}\,|v_{0}(yK)|\,
			\int_{G\smallsetminus{K(\exp\widetilde{\Omega}_{t})K}}\diff{g}\,
			\widetilde{\delta}(g)^{-\frac12}\,Q_{t}^{\sigma, 0}(y^{-1}g)\\[5pt]
			&\lesssim\,
			\underbrace{\vphantom{
					\int_{S\smallsetminus{K(\exp\widetilde{\Omega}_{t}')K}}}
				\int_{S}\textrm{d}_{r}{y}\,|\widetilde{v}_{0}(y)|
			}_{\|\widetilde{v}_{0}\|_{L^{1}(S)}}\,
			\underbrace{
				\int_{S\smallsetminus{K(\exp\widetilde{\Omega}_{t}'')K}}
				\textrm{d}_{r}{g}\,\widetilde{Q}_{t}^{\sigma}(g)
			}_{\longrightarrow\,0}.\
		\end{align*}
		This concludes the proof of the extension problem asymptotics in $L^{1}$ for the 
		distinguished Laplacian $\widetilde{\Delta}$ on $S$ and for initial data
		$\widetilde{v}_{0}\in\mathcal{C}_{c}(S)$.
	\end{proof}
	
	\subsection{Heat asymptotics in $L^{\infty}$
		for compactly supported initial data}
	
	We first recall the following lemma, which allows us to compare the middle	components occurring in the Iwasawa decomposition and in the Cartan decomposition.
	
	\begin{lemma}\cite[Lemma 4.8]{APZ2023}
		For all $g\in{G}$, we have
		\begin{align}\label{S4 ineq Iwasawa Cartan}
			\langle{\rho,A(g)}\rangle\,
			\le\,\langle{\rho,g^{+}}\rangle
		\end{align}
		where $A(g)$ denotes the $\mathfrak{a}$-component of $g$ in the Iwasawa 
		decomposition and $g^{+}$ denotes its $\overline{\mathfrak{a}^{+}}$-component
		in the Cartan decomposition.
	\end{lemma}

	In the following two propositions we collect some elementary properties of the extension problem kernel. The first one clarifies the lower
	and the upper bounds of $\widetilde{Q}_{t}^{\sigma}$, while the second one describes
	its critical region for the $L^{\infty}$ norm.
	
	\begin{proposition}\label{S4 htilde estimate proposition}
		The kernel $\widetilde{Q}_{t}^{\sigma}$ associated with the extension problem for the distinguished Laplacian
		satisfies
		\begin{align}\label{S4 htilde estimate}
			\|\widetilde{Q}_{t}^{\sigma}\|_{L^{\infty}(S)}\,
			\asymp\,t^{-\ell-|\Sigma_{r}^{+}|}
		\end{align}
		for $t$ large enough.
	\end{proposition}

	\begin{proof}
		Using the global estimates \eqref{S2 global estimate phi0} and 
		\eqref{kernel bounds S}, we have
		\begin{align}\label{S4 htilde1}
			\widetilde{Q}_{t}^{\sigma}(g)\,
			\asymp\,
			t^{2\sigma}\,	e^{-\langle{\rho,A(g)}\rangle}e^{-\langle{\rho,g^{+}}\rangle}\,(t+|g|)^{-\ell-2|\Sigma_r^+|-2\sigma}
			\Big\lbrace{
				\prod\nolimits_{\alpha\in\Sigma_{r}^{+}}
				1+ \langle \alpha,g^{+}\rangle 
			}\Big\rbrace
		\end{align}
		We obtain first the lower bound in \eqref{S4 htilde estimate}
		by evaluating the right hand side of \eqref{S4 htilde1} at
		$g_{0}=\exp(-t\rho)$ and by observing that
		\begin{align*}
			A(g_{0})\,=\,-t\rho
			\qquad\textnormal{and}\qquad
			g_{0}^{+}\,=\,t\rho.
		\end{align*}
		For the upper bound, notice that
		\begin{align}\label{S4 htilde2}
			e^{-\langle{\rho,A(g)}\rangle}e^{-\langle{\rho,g^{+}}\rangle}\,
			\le\,1
		\end{align}
		according to \eqref{S4 ineq Iwasawa Cartan}, and that
		\begin{align*}
			(t+|g|)^{-|\Sigma_r^+|}
			\Big\lbrace{
				\prod\nolimits_{\alpha\in\Sigma_{r}^{+}}
				1+ \langle \alpha,g^{+}\rangle 
			}\Big\rbrace\lesssim 1, \qquad t^{2\sigma}\,(t+|g|)^{-\ell-|\Sigma_r^+|-2\sigma}\lesssim t^{-\ell-|\Sigma_r^+|}
		\end{align*}
		for $t$ large enough, whence the claim follows from \eqref{S4 htilde1}.
	\end{proof}

	\begin{proposition}
		The fractional Poisson kernel $\widetilde{Q}_{t}^{\sigma}$ concentrates asymptotically
		in the same critical region for the $L^{\infty}$ norm as for the $L^1$ norm.
		In other words,
		\begin{align*}
			t^{\ell+|\Sigma_{r}^{+}|}\,
			\| \widetilde{Q}_{t}^{\sigma}\|_{L^{\infty}
				(S\smallsetminus{K(\exp\widetilde{\Omega}_{t})K})}\,
			\longrightarrow\,0
			\qquad\textnormal{as}\,\,\,t\rightarrow +\infty.
		\end{align*}
	\end{proposition}
	
	\begin{proof}
		Let us study the sup norm of $\widetilde{Q}_{t}^{\sigma}$ outside the critical region. 
		From  \eqref{S4 htilde1} and \eqref{S4 htilde2} we deduce that 
		\begin{align}\label{S4 htilde1'}
			t^{\ell+|\Sigma_{r}^{+}|}\,
			\widetilde{Q}_{t}^{\sigma}(g)\,
			\lesssim\,
			t^{2\sigma+\ell+|\Sigma_{r}^{+}|}(t+|g^{+}|)^{-2\sigma-\ell-2|\Sigma_{r}^{+}|}(1+|g^{+}|^{|\Sigma_{r}^{+}|}).
		\end{align}
		\noindent\textit{Case 1}: Assume that $|g^{+}|<t^{1-\varepsilon}$. Then $t+|g^{+}|\asymp t$ and $(1+|g^{+}|)^{|\Sigma_{r}^{+}|}\lesssim t^{(1-\varepsilon)|\Sigma_{r}^{+}|}$.
		Thus we deduce from \eqref{S4 htilde1'} that
		\begin{align*}
			t^{\ell+|\Sigma_{r}^{+}|}\,
			\widetilde{Q}_{t}^{\sigma}(g)\,
			\lesssim\,
			t^{-\varepsilon|\Sigma_{r}^{+}|}
		\end{align*}
		which tends to $0$.
		
		\noindent\textit{Case 2}: Assume that $|g^{+}|>t^{1+\varepsilon}$.
		Then $t+|g^{+}|\asymp |g^{+}|$, therefore
		\begin{align*}
			t^{\ell+|\Sigma_{r}^{+}|}\,
			\widetilde{Q}_{t}^{\sigma}(g)\,
			\lesssim\,
			t^{2\sigma+\ell+|\Sigma_{r}^{+}|}|g^{+}|^{-2\sigma-\ell-|\Sigma_{r}^{+}|}\lesssim t^{-\varepsilon(\ell+|\Sigma_{r}^{+}|+2\sigma)}
		\end{align*}
		which tends to $0$. This completes the proof.
	\end{proof}
	
	\vspace{5pt}
	Finally, let us prove the remaining part of \cref{S1 Main thm 2}.
	\begin{proof}[Proof of \eqref{S1 Linf disting} in \cref{S1 Main thm 2}.]
		Fix $0<\varepsilon <\frac{1}{\nu+2\sigma-1}$. Consider the function $$t\mapsto \epsilon(t), \quad \epsilon(t)=t^{-1+\varepsilon(\nu+2\sigma-1)}\longrightarrow 0, \quad \text{as} \quad t\rightarrow +\infty.$$	In the critical region $S\cap{K(\exp\widetilde{\Omega}_{t})K}$, we have
		\begin{align*}
			|\widetilde{v}(t,g)-\widetilde{M}(g)\widetilde{Q}_{t}^{\sigma}(g)|\,
			\le\,\widetilde{Q}_{t}^{\sigma}(g)\,
			\int_{|y|<\,\xi}\diff{g}\,|v_{0}(yK)|\,
			\Big|{
				\frac{Q_{t}^{\sigma,0}(y^{-1}g)}{Q_{t}^{\sigma,0}(g)}
				-\frac{\varphi_{0}(y^{-1}g)}{\varphi_{0}(g)}
			}\Big|
		\end{align*}
		with 
		\begin{align*}
			\Big|{
				\frac{Q_{t}^{\sigma, 0}(y^{-1}g)}{Q_{t}^{\sigma, 0}(g)}
				-\frac{\varphi_{0}(y^{-1}g)}{\varphi_{0}(g)}
			}\Big|\,
			\lesssim\,\epsilon(t)
		\end{align*}
		according to \eqref{S4 difference} 
		and to \cref{S4 Lemma ratios difference}.
		Then we deduce from \eqref{S4 htilde estimate} that
		\begin{align*}
			t^{\ell+|\Sigma_{r}^{+}|}\,
			|\widetilde{v}(t,g)-\widetilde{M}(g)\widetilde{Q}_{t}^{\sigma}(g)|\,
			\lesssim\,\epsilon(t)
			\qquad\forall\,g\in{S\cap{K(\exp\widetilde{\Omega}_{t})K}}
		\end{align*}
		where the right-hand side tends to $0$ as $t\rightarrow +\infty$.
		Outside the critical region, we estimate separately $\widetilde{v}(t,g)$
		and $\widetilde{M}(g)\widetilde{Q}_{t}^{\sigma}(g)$. On the one hand, we know that
		$\widetilde{M}(g)$ is a bounded function and that 
		$\widetilde{Q}_{t}^{\sigma}(g)=\mathrm{o}(t^{-\ell-|\Sigma_{r}^{+}|})$. Then 
		$t^{\ell+|\Sigma_{r}^{+}|}\widetilde{M}(g)\widetilde{Q}_{t}^{\sigma}(g)$ 
		tends to $0$ as $t\rightarrow+\infty$.
		On the other hand, since $g\notin{K(\exp\widetilde{\Omega}_{t})K}$ and 
		$|y|<\xi$ imply that $g^{-1}y\notin{K(\exp\widetilde{\Omega}_{t}'')K}$
		(see \eqref{S4 Omega''}), we obtain
		\begin{align*}
			|\widetilde{v}(t,g)|\,
			\lesssim\,\int_{G}\diff{y}\,
			|\widetilde{v}_{0}(yK)|\,|\widetilde{Q}_{t}^{\sigma}(g^{-1}y)|\,
		\end{align*}
		which is $\textrm{o}(t^{-\ell-|\Sigma_{r}^{+}|})$ outside the critical region.
		In conclusion,
		\begin{align*}
			t^{\ell+|\Sigma_{r}^{+}|}
			\|\widetilde{v}(t,\,\cdot\,)-
			\widetilde{M}\,\widetilde{Q}_{t}^{\sigma}\|_{L^{\infty}(S)}\,
			\longrightarrow\,0
		\end{align*}
		as $t\rightarrow +\infty$.
	\end{proof}
	
	The result for the $L^{p}$ norm follows by convexity.
	\begin{corollary}
		The solution $\widetilde{v}$ to the Cauchy problem \eqref{HE S} with initial 
		data $\widetilde{v}_{0}\in\mathcal{C}_{c}(S)$ satisfies
		\begin{align}\label{S4 Lp disting}
			t^{\frac{\ell+|\Sigma_{r}^{+}|}{p'}}
			\|\widetilde{v}(t,\,\cdot\,)-
			\widetilde{M}\,\widetilde{Q}_{t}^{\sigma}\|_{L^{p}(S)}\,
			\longrightarrow\,0
			\qquad\textnormal{as}\quad\,t\rightarrow +\infty,
		\end{align}
		for all $1<p<\infty$.
	\end{corollary}
	
	\subsection{Asymptotics for other initial data}\label{Subsect other data}
	We have obtained above the long-time asymptotic convergence in $L^{p}$
	($1\le{p}\le\infty$) for the extension problem 
	with compactly supported initial data.
	The following corollaries give some other functional spaces for which the convergence is true, but the question regarding the full $L^1(S)$ class remains open.
	
	\begin{corollary}
		The asymptotic convergences \eqref{S1 L1 disting} and \eqref{S1 Linf disting}, 
		hence \eqref{S4 Lp disting}, still hold with initial data
		$\widetilde{v}_{0}=\widetilde{\delta}^{\tfrac12}v_{0}\in{L}^{1}(S)$ when 
		$v_{0}$ is bi-$K$-invariant.
	\end{corollary}

	\begin{corollary}
		The asymptotic convergences \eqref{S1 L1 disting} and \eqref{S1 Linf disting}, 
		hence \eqref{S4 Lp disting}, 
		still hold with no bi-$K$-invariance condition but under the assumption
		\begin{align}\label{S4 other data assumption}
			\int_{G}\diff{g}\,|v_{0}(gK)|e^{\langle{\rho,g^{+}}\rangle}\,
			<\,\infty.
		\end{align}
	\end{corollary}
	
	The proofs of the above corollaries are similar to those of \cite[Corollary 4.12]{APZ2023} and \cite[Corollary 4.13]{APZ2023}, respectively, thus omitted.

	
	\vspace{10pt}\noindent\textbf{Acknowledgments.}
	
	The author would like to thank M. Bhowmik and S. Pusti for interesting discussions. This work is funded by the Deutsche Forschungsgemeinschaft (DFG, German Research Foundation)--SFB-Gesch{\"a}ftszeichen --Projektnummer SFB-TRR 358/1 2023 --491392403.  This work was also partially supported by the Hellenic Foundation for Research and Innovation, Project HFRI-FM17-1733.  
	
	\bigskip
	
	

	\printbibliography
	

	\vspace{10pt}
	\address{
		\noindent\textsc{Effie Papageorgiou:}
		\href{mailto:papageoeffie@gmail.com}
		{papageoeffie@gmail.com}\\
		Institut f{\"u}r Mathematik, Universit\"at Paderborn, Warburger Str. 100, D-33098
		Paderborn, Germany}

\end{document}

%% file: Concentration1Lnew.tex

\begin{tikzpicture}[line cap=round,line join=round,>=triangle 45,x=1cm,y=1cm,scale=2.5]
	\clip(-1,-1.25) rectangle (5,1.25);
	\draw [line width=0.0pt,color=blue,fill=blue,fill opacity=0.31] (0,0) circle (3cm);
	\fill[line width=0pt,color=white,fill=white,fill opacity=1] (0,0) -- (4.0654939275664805,0.8360139301518585) -- (0,4) -- (-4,1) -- (-4,-1) -- (0,-4) -- (4.183307817025808,-0.8586935566861592) -- cycle;
	
	\draw [line width=0.5pt] (0,0) circle (2cm);
	\draw [line width=0.5pt] (0,0) circle (3cm); 
	\fill[line width=0pt,color=white,fill=white,fill opacity=1] (0,0) -- (4,1.5) -- (0,4) -- (-4,1) -- (-4,-1) -- (0,-4) -- (4,-1.5) -- cycle;
	
	\draw [line width=0pt,color=white,fill=white,fill opacity=1] (0,0) circle (1.99cm);
	
	\draw [line width=0.5pt] (0,0)-- (4.0654939275664805,0.8360139301518585);
	\draw [line width=0.5pt] (0,0)-- (4.0654939275664805,-0.8360139301518585); 
	
\draw [line width=0.5pt] (0,0)-- (4.0654939275664805,1.1);
\draw [line width=0.5pt] (0,0)-- (4.0654939275664805,-1.1);  
	
	\draw [line width=0.5pt] (0,0)-- (4,0);
	\draw [color=white,fill=white,fill opacity=1] (1.5,-1) rectangle (2,-1);
	\draw [color=white,fill=white,fill opacity=1] (1.5,1) rectangle (2,1);
	\draw [color=white,fill=white,fill opacity=1] (2.5,-1) rectangle (3,-1.5);
	\draw [color=white,fill=white,fill opacity=1] (2.5,1) rectangle (3,1.5);

	\draw [line width=0.5pt] (0.75,0) arc (0:30:0.29);
	\draw [line width=0.5pt] (1.1,0) arc (0:30:0.55);
	
	\draw [dashed, line width=0.5pt] (0,0)-- (4.6,-1.5);
	\draw [dashed, line width=0.5pt] (0,0)-- (4.6,1.5);
	\draw[decoration={brace,mirror,raise=5pt},decorate]
	(4.2,-1.25) -- node[right=6pt] {$\,\,\overline{\mathfrak{a}^{+}}$} (4.2,1.5);
	
	\begin{scriptsize}
		\draw (2.5,0.1) node{$\Omega_{t}$};
		\draw (1.75,-0.8) node{$|H|=t^{2-\varepsilon}$};
		\draw (2.75,-1.1) node{$|H|=t^{2+\varepsilon}$};
		\draw (4,-0.1) node{$\rho$-axis};
		\draw (0.85,0.09) node{$\gamma(t)$};
		\draw (1.16,0.10) node{$\gamma_0$};
		\draw (4,1.21) node{wall};
	\end{scriptsize}
\end{tikzpicture}